\documentclass[10pt,A4]{article}

\usepackage{amsmath}
\usepackage{amsfonts,amssymb}
\usepackage{A4}
\usepackage[francais]{babel}
\usepackage{epsfig,epsf}

\newcounter{fig}
\def\figcaption #1
  {\addtocounter{fig}{1}
   \begin{center}
   {\sc Fig.} \arabic{fig}. {\it #1} 
    \end{center}
    \addtocontents{lof}{\protect \contentsline {figure}{\protect \numberline {\thefig}{\protect \ignorespaces \protect \pit #1 }}{\thepage}}
   }
\def\og{\leavevmode\raise.24ex\hbox{$\scriptscriptstyle\langle\!\langle\>$}}    % guillemets ouvrants
  \def\fg{\leavevmode\raise.24ex\hbox{$\scriptscriptstyle\>\rangle\!\rangle$}}    % guillemets fermants

\newtheorem{prop}{Proposition}

\newcommand{\fhi}{\varphi}
\newcommand{\ioe}{\leqslant}
\newcommand{\soe}{\geqslant}
\newcommand{\vers}{\rightarrow}
\newcommand{\dist}{{\rm dist}}
\newcommand{\Ssi}{\Longleftrightarrow}
\newcommand{\Ecos}{{\rm Ecos}}
\newcommand{\Esin}{{\rm Esin}}

\newcommand{\Hil}{{\mathcal H}}
\newcommand{\Dbv}{{\mathcal W}}

\newcommand{\Nat}{{\mathbb N}}
\newcommand{\Int}{{\mathbb Z}}

\newcommand{\Real}{{\mathbb R}}
\newcommand{\Com}{{\mathbb C}}

\newcommand{\Res}{{\rm Res}}

\newcommand{\fin}{\hfill$\Box$}
\newcommand{\dem}{\noindent {\bf D\'emonstration\ }}
\newcommand{\fine}{\tag*{\mbox{$\Box$}}}

\title{Sur l'autocorr\'elation multiplicative de la fonction \\ {\og
partie fractionnaire \fg}}

\author{Luis B\'AEZ-DUARTE, Michel BALAZARD, Bernard LANDREAU et Eric SAIAS}
\date{}
\begin{document}
\maketitle

Le texte qui suit est un document de travail, contenant les
pr\'erequis et les d\'etails des calculs menant aux r\'esultats de
notre article, 

\begin{center}
{\it \'Etude de l'autocorr\'elation multiplicative de la fonction {\og
partie fractionnaire \fg}}
\end{center}
\`a para\^{\i}tre au Ramanujan Journal.

\newpage
 
\tableofcontents

\newpage

\section{Transformation de Mellin}\label{t94}

\subsection{Convergence absolue}

Soit $]a,b[$ un intervalle ouvert non vide de $\Real \,$, et $\Dbv (a,b)$
l'espace vectoriel des fonctions $f : \, ]0,+\infty[ \vers \Com$,
mesurables au sens de Lebesgue, et telles que, pour tout $\sigma \in ]a,b[$,
$$
\int_0^{+\infty} |f(x)|x^{\sigma-1} dx < +\infty.
$$
Observons que $\Dbv (a,b) \subset \Dbv (a',b')$ si $]a',b'[ \subset
]a,b[$.

Si $f \in \Dbv (a,b)$, la {\bf transform\'ee de Mellin} de $f$ est la fonction $Mf$ d\'efinie par
$$
Mf(s) := \int_0^{+\infty} f(x)x^{s-1} dx.
$$
Elle est d\'efinie et holomorphe dans la bande verticale $a < \Re s < b$
(qui peut \^etre un demi-plan, ou m\^eme le plan tout entier). Elle est
born\'ee dans toute bande verticale $\alpha \ioe \Re s \ioe
\beta$, o\`u $a < \alpha  \ioe \beta < b$.

\begin{prop}[Injectivit\'e de la transformation de Mellin]
Si $f \in \Dbv (a,b)$ et $Mf=0$, alors $f=0$ p.p.
\end{prop}

\smallskip\noindent 

{\bf Exemples}

$\bullet$  $f(x)=e^{-x}$;  $f \in \Dbv (0, +\infty)$ et
$Mf(s)=\Gamma(s)$.
\smallskip

$\bullet$  $f(x)=\sin x$;  $f \in \Dbv (-1,0)$ et
$Mf(s)=\sin \frac{\pi s}{2}\Gamma(s)$.
\smallskip

$\bullet$  $f(x)=\cos x$ -1;  $f \in \Dbv (-2,0)$ et
$Mf(s)=\cos \frac{\pi s}{2}\Gamma(s)$.

$\bullet$  Posons
$$
f(x)= \begin{cases}
x^{-\rho} (\log x) ^k&  0<x \ioe 1,\\ 
0 & x > 1,
\end{cases} 
$$
o\`u $\rho \in \Com$ et $k \in \Nat$. On a : $f \in \Dbv (\Re
\rho, + \infty)$ et
$$
Mf(s)= \frac{(-1)^k k!}{(s-\rho)^{k+1}}.
$$
On constate que $Mf$ est une fraction rationnelle et que $f(x)= \Res
\bigl ( Mf(s)x^{-s},\rho \bigr )$ pour $0<x \ioe 1$.
\smallskip

$\bullet$  Posons
$$
f(x)= \begin{cases}
0 & 0<x \ioe 1,\\ 
x^{-\rho} (\log x) ^k&   x > 1,
\end{cases} 
$$
o\`u $\rho \in \Com$ et $k \in \Nat$. On a : $f \in \Dbv (-\infty,\Re
\rho)$ et
$$
Mf(s)= -\frac{(-1)^k k!}{(s-\rho)^{k+1}}.
$$
On constate que $Mf$ est l'oppos\'ee de la fraction rationnelle
pr\'ec\'edente et que $f(x)= -\Res
\bigl ( Mf(s)x^{-s},\rho \bigr )$ pour $x > 1$
\smallskip

Appelons polyn\^ome g\'en\'eralis\'e toute somme finie
\begin{equation}\label{t79}
\sum c x^{-\rho} (\log x) ^k,
\end{equation}
o\`u les $c$ et les $\rho$ sont des nombres complexes, et les $k$ des
nombres entiers naturels. Associons \`a ce polyn\^ome g\'en\'eralis\'e
la fraction rationnelle
$$
F(s):= \sum c \frac{(-1)^k k!}{(s-\rho)^{k+1}}.
$$
C'est une fraction rationnelle nulle \`a l'infini. Inversement, toute
fraction rationnelle nulle \`a l'infini s'\'ecrit sous cette forme et
on peut lui associer le polyn\^ome g\'en\'eralis\'e (\ref{t79}). La
correspondance ainsi d\'efinie est bijective et les deux derniers
exemples conduisent \`a la proposition suivante.
\begin{prop}\label{t78}
Soit 
$$
P(x)=\sum c x^{\rho} (\log x) ^k
$$
un polyn\^ome g\'en\'eralis\'e, et
$$
F(s):= \sum c \frac{(-1)^k k!}{(s-\rho)^{k+1}}
$$
la fraction rationnelle associ\'ee. On pose $a= \min \Re \rho$ et $b=
\max \Re \rho$.

On a, pour tout $x>0$,
$$
P(x)= \sum_{\rho} \Res \bigl ( F(s) x^{-s}, \rho \bigr ),
$$
o\`u la somme porte sur les p\^oles de $F$.

D'autre part, la fonction
$$
f(x)= \begin{cases}
 P(x) \quad &0 < x \ioe 1\\
  0   \quad  &x>1,
      \end{cases}
$$
appartient \`a $\Dbv(b,+\infty)$ et v\'erifie $Mf=F$ ; et la fonction
$$
g(x)= \begin{cases}
 0 \quad &0 < x \ioe 1\\
 P(x)  \quad  &x>1,
      \end{cases}
$$
appartient \`a $\Dbv(-\infty,a)$ et v\'erifie $Mg=-F$.
\end{prop}

\begin{prop}
Si $f \in \Dbv (a,b)$ et $\rho \in \Com$, alors la fonction $g$
d\'efinie par $g(t) := t^{\rho}f(t)$ appartient \`a \\
$\Dbv (a-\Re
\rho,b- \Re \rho)$ et
$$
Mg(s) = Mf (s+\rho).
$$
\end{prop}

\begin{prop}
Si $f \in \Dbv (a,b)$, alors la fonction $g$
d\'efinie par $g(t) := f(1/t)$ appartient \`a $\Dbv (-b,-a)$ et
$$
Mg(s)=Mf(-s).
$$
\end{prop}

\begin{prop}
Si $f \in \Dbv (a,b)$, alors $\overline{f} \in \Dbv (a,b)$ et
$$
M\overline{f}(s) = (Mf)^*(s) := \overline{Mf(\overline{s})}.
$$
\end{prop}

\begin{prop}
Si $f \in \Dbv (a,b)$ et $\lambda > 0$, alors la fonction $g$
d\'efinie par $g(t)=f(\lambda t)$ appartient \`a $\Dbv (a,b)$ et
$$
Mg(s)=\lambda^{-s} Mf(s).
$$
\end{prop}

\begin{prop}
Si $f$ et $g$ appartiennent \`a $\Dbv (a,b)$, la fonction $f*g$
d\'efinie presque partout par
$$
f*g(x) := \int_0^{+\infty} f(t) g(x/t) \frac{dt}{t},
$$
appartient \`a $\Dbv (a,b)$ et v\'erifie
$$
M(f*g) = Mf \cdot Mg.
$$
\end{prop}

Observons que $f*g$ est continue si $f$ ou $g$ est continue et \`a
support compact dans $]0,+\infty[$.
\smallskip

{\bf Exemples}

$\bullet$ Si $f \in \Dbv (a,b)$ et
$$
g(x)=\begin{cases}
x^{-\rho} &  0<x \ioe 1,\\ 
0 & x > 1,
\end{cases} 
$$
o\`u $\Re \rho \ioe a$, alors
\begin{align*}
f*g(x) &= \int_x^{+\infty} f(t) \left ( \frac{x}{t} \right ) ^{-\rho}
\frac{dt}{t} \\
&= x^{-\rho} \int_x^{+\infty} f(t) t^{\rho -1}dt.
\end{align*}
On a $f*g \in \Dbv (a,b)$ et
$$
M(f*g)(s) = \frac{Mf(s)}{s-\rho}.
$$
\smallskip

$\bullet$ Si $f \in \Dbv (a,b)$ et
$$
g(x)=\begin{cases}
0 &0<x \ioe 1,\\ 
x^{-\rho} &  x > 1,
\end{cases} 
$$
o\`u $\Re \rho \soe b$, alors
\begin{align*}
f*g(x) &= \int_0^x f(t) \left ( \frac{x}{t} \right ) ^{-\rho}
\frac{dt}{t} \\
&= x^{-\rho} \int_0^x f(t) t^{\rho -1}dt.
\end{align*}
On a $f*g \in \Dbv (a,b)$ et
$$
M(f*g)(s) = -\frac{Mf(s)}{s-\rho}.
$$
\smallskip

\begin{prop}\label{t93}
Si $f \in \Dbv (a,b)$ et si les nombres complexes $c_k$ et les nombres
r\'eels positifs $\lambda_k$ sont tels que
$$
\sum |c_k| \lambda_k^{-\sigma} < +\infty
$$
pour tout $\sigma \in ]a,b[$, alors la fonction $g$ d\'efinie presque
partout par
$$
g(t) = \sum c_k f(\lambda_k t)
$$
appartient \`a $\Dbv (a,b)$ et v\'erifie
$$
Mg(s) = Mf(s) \cdot \sum c_k \lambda_k^{-s}.
$$
\end{prop}

\begin{prop}[Inversion de Mellin]
Soit $f \in \Dbv(a,b)$. L'ensemble des $c \in ]a,b[$ tels que $\tau
\mapsto Mf(c+i \tau)$ est dans $L^1(\Real)$ est un intervalle
(\'eventuellement vide). Si $c$ est choisi dans cet ensemble, on a
$$
f(t) = \frac{1}{2 \pi i} \int_{\sigma=c} Mf(s) t^{-s} ds \quad p.p.
$$
\end{prop}

Nous dirons d'une fonction $F(s)$ m\'eromorphe dans la bande verticale
$a < \Re s < b$ qu'elle est {\bf \`a croissance polyn\^omiale} si elle
n'a qu'un nombre fini de p\^oles et si, pour tous $\alpha$, $\beta$
tels que $a < \alpha \ioe \beta < b$, existent deux nombres r\'eels
$T_0 >0$ et $K$, pouvant d\'ependre de $\alpha$ et $\beta$, tels que
$$
|F(\sigma + i \tau)| \ioe |\tau|^K, \quad \alpha \ioe \sigma \ioe
 \beta, \, |\tau| \soe T_0.
$$ 

\begin{prop}\label{t97}
Soit $]a,b[$ et $]a',b'[$ deux intervalles ouverts non vides, avec $b
\ioe a'$. Soit $F(s)$ une fonction m\'eromorphe et \`a croissance
polyn\^omiale dans la bande verticale $a < \Re s < b'$. Soit enfin $f
\in \Dbv(a,b)$ et $g \in \Dbv(a',b')$ deux fonctions telles que
$$
Mf(s) = F(s), \quad a < \Re s < b;
$$
$$
Mg(s) = F(s), \quad a' < \Re s < b'.
$$

Alors $g(t)-f(t)$ coincide presque partout avec le polyn\^ome
g\'en\'eralis\'e
$$
P(t) = \sum_{\rho} \Res \bigl ( F(s)t^{-s}, \rho \bigr ),
$$
o\`u la
somme porte sur les p\^oles de $F$.
\end{prop}
{\bf D\'emonstration}

Soit $A(s)$ la fraction rationnelle obtenue en ajoutant les parties
polaires de $F(s)$ en tous ses p\^oles. On a :
$$
\Res \bigl ( F(s)t^{-s}, \rho \bigr )=\Res \bigl ( A(s)t^{-s}, \rho
\bigr )
$$
pour tout $t>0$ et tout p\^ole $\rho$ de $F$.

D'apr\`es la proposition \ref{t78}, la fonction $g - P \chi_{(0,1)}$
appartient \`a $\Dbv (a',b')$ et admet $F-A$ comme transform\'ee de
Mellin ; et la fonction $f +P\chi_{(1,+\infty)}$ appartient \`a $\Dbv
(a,b)$ et admet $F-A$ comme transform\'ee de Mellin. D'autre part,
$F-A$ est holomorphe et \`a croissance polyn\^omiale dans la bande $a
< \Re s < b'$. Si le r\'esultat est connu dans le cas o\`u $F$ est
holomorphe, le cas g\'en\'eral en d\'ecoule donc.

Supposons donc $F$ holomorphe. Comme $F$ est une transform\'ee de
Mellin dans les deux bandes verticales $a
< \Re s < b$ et $a'
< \Re s < b'$, le principe de Phragmen-Lindel\"of prouve que $F$ est
born\'ee dans toute bande verticale $\alpha \ioe \Re s \ioe \beta$,
o\`u $a < \alpha \ioe \beta <b'$. Soit $\fhi$ une fonction deux fois
contin\^ument d\'erivable et \`a support compact dans $]0,+\infty[$ ; les fonctions
$\fhi*f$ et $\fhi*g$ appartiennent respectivement \`a $\Dbv (a,b)$ et
$\Dbv (a',b')$ et ont pour transform\'ees de Mellin $F(s) \cdot M \fhi
(s)$ dans les deux cas. Comme $M \fhi (s)$ est $O \bigl ( (1 +|\tau
|)^{-2} \bigr )$ sur toute droite verticale (et m\^eme uniform\'ement
dans toute bande verticale de largeur finie), la formule d'inversion
de Mellin donne pour presque tout $t$
$$
(f*\fhi) (t) = \frac{1}{2 \pi i} \int_{\sigma=c} M \fhi (s) \cdot F(s)
t^{-s} ds, \quad a<c<b,
$$
et
$$
(g*\fhi) (t) = \frac{1}{2 \pi i} \int_{\sigma=c'} M \fhi (s) \cdot F(s)
t^{-s} ds, \quad a'<c'<b'.
$$

Le th\'eor\`eme de Cauchy prouve alors que $(f*\fhi)(t)=(g*\fhi)(t)$
presque partout, et m\^eme pour tout $t$, par continuit\'e. Comme
cette \'egalit\'e est valable pour toute $\fhi$, on a $f(t)=g(t)$ p.p.\fin

\subsection{Transformation de Mellin-Plancherel}\label{t113}

Soit $f \in L^2(0,+\infty)$. Pour $T>1$, posons :
$$
F_T(s) := \int_{1/T}^T f(t)t^{s-1}dt.
$$

La restriction de $F_T$ \`a la droite $\Re s =1/2$ appartient \`a
$L^2(1/2 + \Real i; d \tau/ 2 \pi)$ et, dans cet espace, $F_T$ a une
limite quand $T$ tend vers l'infini. Nous noterons encore $Mf(s)$
cette fonction, dite {\bf transform\'ee de Mellin-Plancherel}
de $f$ (ou simplement : transform\'ee de Mellin), d\'efinie presque partout sur la droite $\Re s = 1/2$. Le
th\'eor\`eme de Plancherel affirme que l'application $f \mapsto Mf$
est un op\'erateur unitaire entre les deux espaces de Hilbert
$L^2(0,+\infty)$ et $L^2(1/2 + \Real i; d \tau/ 2 \pi)$. 

Si $\alpha \in \Real$, on d\'efinit plus g\'en\'eralement la
transformation de Mellin-Plancherel sur l'espace \\
$L^2(0,+\infty;t^{\alpha}dt)$, \`a valeurs dans $L^2 \bigl ( (\alpha +1)/2 +
\Real i; d \tau/ 2 \pi \bigr )$. C'est encore un op\'erateur unitaire.

\newpage

\section{La fonction \og partie
fractionnaire \fg }

\subsection{Une formule sommatoire}

La proposition suivante\footnote{Cf. Lemma 7 dans : {\sc G. H. Hardy et
J. E. Littlewood}, {\it Some problems of Diophantine approximation :
The lattice-points of a right-angled triangle (Second memoir)}, Abh. Math. Sem. Hamb. Un. {\bf 1} (1922), 212-249.}
illustre un principe g\'en\'eral de sym\'etrie.

\begin{prop}\label{t85}
Soit $\theta$ et $x$ deux nombres r\'eels positifs, $f(n)$ une
fonction \`a valeurs complexes, d\'efinie pour $n$ entier, $0 \ioe n
\ioe \theta x$, et $g(m)$ une
fonction \`a valeurs complexes, d\'efinie pour $m$ entier, $0 \ioe m
\ioe x$. On suppose de plus que $f(0)=g(0)=0$. 

Si $\theta$ est irrationnel, on a
$$
\sum_{ 1 \ioe m \ioe x} f \bigl ( \left \lfloor m \theta \right \rfloor
\bigr ) \bigl ( g(m) -g(m-1) \bigr ) + \sum_{ 1 \ioe n \ioe \theta x} g \bigl ( \left \lfloor n/ \theta \right \rfloor
\bigr ) \bigl ( f(n) - f(n-1) \bigr ) = f \bigl ( \left \lfloor
\theta x \right \rfloor
\bigr ) g \bigl ( \left \lfloor x \right \rfloor
\bigr ).
$$

Si $\theta = p/q$, o\`u $p$ et $q$ sont des
nombres entiers positifs premiers entre eux, ce r\'esultat est \`a remplacer par
$$
 f \bigl ( \left \lfloor
\theta x \right \rfloor
\bigr ) g \bigl ( \left \lfloor x \right \rfloor
\bigr ) + \sum_{1 \ioe k \ioe x/q} \bigl ( g(kq) - g(kq-1) \bigr ) \bigl
( f(kp)-f(kp-1) \bigr ).
$$
\end{prop}
{\bf D\'emonstration}

Voyons d'abord comment le r\'esultat pour $p/q$ d\'ecoule du cas o\`u
$\theta$ est irrationnel. Faisons tendre $\theta$ vers $p/q$ par
valeurs sup\'erieures et irrationnelles. La limite de
$$
\sum_{ 1 \ioe m \ioe x} f \bigl ( \left \lfloor m \theta \right \rfloor
\bigr ) \bigl ( g(m) -g(m-1) \bigr ) -f \bigl ( \left \lfloor
\theta x \right \rfloor
\bigr ) g \bigl ( \left \lfloor x \right \rfloor
\bigr )
$$
s'obtient en rempla\c{c}ant $\theta$ par $p/q$, alors que celle de
$$
\sum_{ 1 \ioe n \ioe \theta x} g \bigl ( \left \lfloor n/ \theta \right \rfloor
\bigr ) \bigl ( f(n) - f(n-1) \bigr ) 
$$
vaut
$$
\sum_{\substack{
        1 \ioe n \ioe px/q\\
        nq/p \not \in  \Nat}}
g \bigl ( \left \lfloor nq/p \right \rfloor
\bigr ) \bigl ( f(n) - f(n-1) \bigr )  + \sum_{\substack{
        1 \ioe n \ioe px/q\\
        nq/p  \in  \Nat}}
g \left ( \frac{nq}{p} -1
\right ) \bigl ( f(n) - f(n-1) \bigr )
$$
$$
= \sum_{
        1 \ioe n \ioe px/q}
g \bigl ( \left \lfloor nq/p \right \rfloor
\bigr ) \bigl ( f(n) - f(n-1) \bigr )  + \sum_{\substack{
        1 \ioe n \ioe px/q\\
        nq/p  \in  \Nat}}
\left ( g \left ( \frac{nq}{p} -1
\right )- g \left ( \frac{nq}{p} \right ) \right ) \bigl ( f(n) -
        f(n-1) \bigr ).
$$
La derni\`ere somme se r\'ecrit en posant $n=kp$ :
$$
\sum_{k \ioe x/q} \bigl ( g(kq-1) - g(kq) \bigr ) \bigl
( f(kp)-f(kp-1) \bigr ),
$$
d'o\`u le r\'esultat.

Si $\theta$ est irrationnel, on a :
$$
\sum_{ 1 \ioe m \ioe x} f \bigl ( \left \lfloor m \theta \right \rfloor
\bigr ) \bigl ( g(m) -g(m-1) \bigr )=  \sum_{ 0 \ioe n \ioe \theta x} f(n) \sum_{\substack{
                                                         1 \ioe m \ioe x\\
                                                         \left \lfloor
m \theta \right \rfloor =n}} \bigl ( g(m) -g(m-1) \bigr ).
$$
Or
\begin{align*}
\left \lfloor
m \theta \right \rfloor =n & \Leftrightarrow n \ioe m \theta < n+1 \\
& \Leftrightarrow \frac{n}{\theta} \ioe m <\frac{n+1}{\theta} \\
& \Leftrightarrow \left \lfloor \frac{n}{\theta} \right \rfloor +1
\ioe m < \left \lfloor \frac{n+1}{\theta} \right \rfloor.
\end{align*}
Si $n \ioe \theta x -1$, cet encadrement entra\^{\i}ne que $m \ioe x$,
alors que si $n= \lfloor \theta x \rfloor$ on a
\begin{align*}
\left \lfloor
m \theta \right \rfloor =n \quad \text{et} \quad m \ioe x \quad  &
\Leftrightarrow \quad 
n \ioe m \theta < n+1 \quad \text{et} \quad m \ioe x \\
& \Leftrightarrow \quad\left \lfloor \frac{n}{\theta} \right \rfloor +1
\ioe m \ioe \lfloor x \rfloor.
\end{align*}
Par cons\'equent on a
$$
\sum_{\substack{
                                                         m \ioe x\\
                                                         \left \lfloor
m \theta \right \rfloor =n}} \bigl ( g(m) -g(m-1) \bigr ) =
                                                         \begin{cases} 
g \left (\left \lfloor \frac{n+1}{\theta} \right \rfloor \right ) - g
                                                         \left (\left
                                                         \lfloor
                                                         \frac{n}{\theta}
                                                         \right
                                                         \rfloor
                                                         \right )&
                                                         \text{si $n \ioe
                                                         \theta x
                                                         -1$},\\
g \bigl ( \lfloor x \rfloor \bigr ) - g\left (\left
                                                         \lfloor
                                                         \frac{n}{\theta}
                                                         \right
                                                         \rfloor
                                                         \right )&
                                                         \text{si $n=
                                                         \lfloor
                                                         \theta x
                                                         \rfloor$}.
							 \end{cases}
$$
Il en r\'esulte que
\begin{align*}
\sum_{ 0 \ioe n \ioe \theta x} f(n) 
\sum_{\substack{
                1 \ioe m \ioe x\\
                \left \lfloor m \theta \right \rfloor =n}} 
\bigl ( g(m) -g(m-1) \bigr ) &= 
\sum_{ 0 \ioe n \ioe \theta x -1}
f(n) \left (g \left (\left \lfloor \frac{n+1}{\theta} \right \rfloor
\right ) 
- g \left (\left \lfloor \frac{n}{\theta} \right \rfloor \right )
\right )\\
& \quad + f \bigl (\lfloor \theta x \rfloor \bigr )
\left ( g \bigl (\lfloor x \rfloor \bigr )
- g \left ( \left \lfloor \frac{\lfloor \theta x \rfloor}{\theta}
\right \rfloor \right ) \right )\\
&= g \bigl (\lfloor x \rfloor \bigr ) f \bigl (\lfloor \theta x \rfloor
                                                         \bigr ) + \sum_{ 1 \ioe n \ioe \theta x} g \bigl ( \left \lfloor n/ \theta \right \rfloor
\bigr ) \bigl ( f(n-1) - f(n) \bigr ),
\end{align*}
comme annonc\'e.\fin

En prenant $f(n)=g(n)=n$ dans la proposition \ref{t85}, on obtient en
particulier le r\'esultat suivant\footnote{{\sc J. J. Sylvester}, {\it
Sur la fonction $E(x)$}, C. R. A. S. {\bf 50} (1860), 732-734.}.

\begin{prop}\label{t86}
Soit $x$ et $\theta$ deux nombres r\'eels positifs. La quantit\'e
$$
\sum_{1 \ioe m \ioe x} \lfloor m \theta \rfloor + \sum_{1 \ioe n \ioe
\theta x} \lfloor n/ \theta \rfloor
$$
vaut $\lfloor x \rfloor  \lfloor \theta x \rfloor$ si $\theta$
est irrationnel, et $\lfloor x \rfloor  \lfloor \theta x \rfloor
+ \lfloor x/q \rfloor$ si $\theta = p/q$, avec $p \in \Nat^*$, $q \in
\Nat^*$, $(p,q)=1$.
\end{prop}

\subsection{La fonction {\og partie fractionnaire \fg}}

La {\bf partie fractionnaire} du nombre r\'eel $x$ est l'unique nombre
r\'eel $u \in [0,1[$ tel que $x-u$ soit entier. On note
$$
\{x\} := u = x - \lfloor x \rfloor,
$$ car l'entier $x-u$ est n\'ecessairement la partie enti\`ere de $x$.

On dispose donc d'une fonction \og partie fractionnaire \fg : $\{ \, \} :
\Real \vers [0,1[$.
\begin{prop}
La fonction {\og partie fractionnaire \fg} est l'unique fonction $f: \Real
\vers \Real$, p\'eriodique de p\'eriode $1$, et v\'erifiant $f(x)=x$
pour $0 \ioe x <1$.
\end{prop}

La fonction {\og partie fractionnaire \fg} poss\`ede donc une s\'erie de
Fourier
$$
\frac{1}{2} - \sum_{n \not = 0} \frac{e^{i 2 \pi n x}}{i 2 \pi n},
$$
qu'on peut r\'ecrire sous la forme
$$
\frac{1}{2} - \sum_{n \soe 1}\frac{\sin 2 \pi n x}{\pi n}.
$$
\begin{prop}\label{t89}
On a
$$ 
\frac{1}{2} - \sum_{n \soe 1}\frac{\sin 2 \pi n x}{\pi n} = 
\begin{cases}
  \{x\}&   \text{si $x \not \in \Int$;}\\
  \frac{1}{2}& \text{si $x \in \Int$.}
\end{cases}
$$

De plus\footnote{{\sc T. H. Gronwall}, {\it \"Uber die Gibbsche
Erscheinung und die trigonometrischen Summen $\sin x + \frac{1}{2}
\sin 2x + \cdots + \frac{1}{n}\sin nx$}, Math. Annalen {\bf 72} (1912), 228-243.}, pour tout nombre entier $N$ on a
$$
0 < \sum_{n=1}^N\frac{\sin 2 \pi n x}{\pi n}  < \frac{1}{\pi}
\int_0^{\pi} \frac{\sin t}{t} dt = 0,589 \dots, \quad 0<x< \frac{1}{2}.
$$
\end{prop}

\begin{prop}\label{t87}
Pour $x \soe 0$, on a
$$
\sum_{n \ioe x} n = \frac{x^2}{2} -x \bigl ( \{ x \} - 1/2 \bigr ) +
\frac{\{ x \}^2}{2} - \frac{\{ x \}}{2}.
$$
\end{prop}
{\bf D\'emonstration}

On a
\begin{align*}
\sum_{n \ioe x} n &= \frac{\lfloor x \rfloor \bigl ( \lfloor x \rfloor
+1 \bigr )}{2}\\
&= \frac{\bigl ( x - \{ x \} \bigr )\bigl ( x - \{ x \} +1 \bigr
)}{2}\\
&= \frac{1}{2} \bigl ( x^2 -2 \{ x \}x +\{ x \}^2 + x - \{ x \} \bigr
)\\
&= \frac{x^2}{2} -x \bigl ( \{ x \} - 1/2 \bigr ) +
\frac{\{ x \}^2}{2} - \frac{\{ x \}}{2}.\fine
\end{align*}

\'Evaluons maintenant une int\'egrale du type de Frullani faisant
intervenir la fonction \og partie fractionnaire \fg.

\begin{prop}\label{t109}
Soit $x$ et $\theta$ deux nombres r\'eels positifs. On a
$$
\int_0^x t^{-2} \bigl ( \{ \theta t \} - \theta \{  t \} \bigr ) dt =
\theta \log \frac{1}{\theta} + \theta \int_x^{\theta x} u^{-2} \{  u
\} du.
$$
\end{prop}
\dem

Posons $\alpha = \min(1,1/\theta)$. On a :
\begin{align*}
\int_0^x t^{-2} \bigl ( \{ \theta t \} - \theta \{  t \} \bigr ) dt &=
\int_{\alpha}^x t^{-2} \bigl ( \{ \theta t \} - \theta \{  t \} \bigr
) dt \\
&= \int_{\alpha}^x t^{-2}\{ \theta t \}dt - \theta \int_{\alpha}^x
t^{-2}\{  t \}dt\\
&= \theta \left ( \int_{\theta \alpha}^{\theta x}u^{-2} \{  u
\} du - \int_{\alpha}^x
u^{-2}\{  u \}du \right )\\ 
&= \theta \left ( \int_{\theta \alpha}^{\alpha} u^{-2} \{  u
\} du + \int_x^{\theta x} u^{-2} \{  u
\} du \right ) \\
&=\theta \log \frac{1}{\theta} + \theta \int_x^{\theta x} u^{-2} \{  u
\} du. \fine
\end{align*}

\subsection{Fonctions de Bernoulli}

On d\'efinit les {\bf polyn\^omes de Bernoulli} $b_n(x)$, $n \soe 1$, par
l'identit\'e formelle
$$
\frac{t e^{xt}}{e^t-1} = \sum_{n=0}^{+\infty}b_n(x) \frac{t^n}{n!}.
$$
On a par exemple
\begin{eqnarray*}
b_0(x) & = & 1 \\
b_1(x) & = & x- \frac{1}{2} \\
b_2(x) & = & x^2 - x+ \frac{1}{6} \\
b_3(x) & = & x^3 - \frac{3}{2}x + \frac{x}{2} \\
b_4(x) & = & x^4 - 2 x^3 + x^2 - \frac{1}{30}.
\end{eqnarray*}

Les {\bf fonctions de Bernoulli} $B_n(x)$ sont les fonctions de variable
r\'eelle, de p\'eriode $1$, d\'efinies par
$$
B_1(x) =\begin{cases}
\{x \} - \frac{1}{2} & \text{si }x \not \in \Int;\\
0 & \text{si }x  \in \Int,
	\end{cases}
$$
et, pour $n \soe 2$,
$$
B_n(x)=b_n(\{x\}).
$$
\begin{prop} 
Pour tout $x \in \Real$, on a
$$
B_1(x) = \frac{\{x\}-\{-x\}}{2}.
$$
\end{prop}

\begin{prop}
$B_n$ a la parit\'e de $n$ :
$$
B_n(-x) = (-1)^n B_n(x).
$$
\end{prop}

\begin{prop}
Pour $n \soe 2$, on a
$$
B_n(x)= \int_0^x nB_{n-1}(t)dt + B_n (0).
$$
Par cons\'equent, $B_n$ est de
classe ${\cal C}^{n-2}$. 
\end{prop}
\begin{prop}
Pour $n \soe 1$, la s\'erie de Fourier de $B_n$ est
$$
-n!\sum_{k \not = 0} \frac{e^{i2\pi k x}}{(i2 \pi k )^n};
$$
elle converge pour tout $x$ vers $B_n(x)$.
\end{prop}

\begin{prop}\label{t88}
Soit $x$ et $\theta$ deux nombres r\'eels positifs. On a
$$
\sum_{1 \ioe n \ioe \theta x} B_1(n/\theta) + \sum_{1 \ioe m \ioe x}
B_1(m \theta) = \frac{1}{2\theta} \bigl ( \{ \theta x \}-\theta\{ x \}
\bigr )^2 + \frac{\theta - 1}{2 \theta} \bigl ( \{ \theta x \}-\theta\{ x \}
\bigr ).
$$
\end{prop}
{\bf D\'emonstration}

Commen\c{c}ons par supposer $\theta$ irrationnel. Nous aurons,
d'apr\`es les propositions \ref{t86} et \ref{t87}, 
\begin{align*}
\sum_{1 \ioe n \ioe \theta x} B_1(n/\theta) + \sum_{1 \ioe m \ioe x}
B_1(m \theta) &= \sum_{1 \ioe n \ioe \theta x} \left ( n/\theta -
\lfloor n/\theta \rfloor - \frac{1}{2} \right ) + \sum_{1 \ioe m \ioe
x} \left ( m\theta -
\lfloor m\theta \rfloor - \frac{1}{2} \right )\\
&= \frac{1}{\theta} \left ( \frac{(\theta x)^2}{2} - \theta x \bigl (
\{ \theta x \} -1/2 \bigr ) +   \frac{ \{ \theta x \}^2}{2} - \frac{
\{ \theta x \}}{2} \right ) + \\
& \quad + \theta \left ( \frac{x^2}{2} -  x \bigl (
\{  x \} -1/2 \bigr ) +   \frac{ \{  x \}^2}{2} - \frac{
\{  x \}}{2} \right ) + \\ 
& \quad - \lfloor x \rfloor \lfloor \theta x \rfloor - \frac{1}{2}
\lfloor \theta x \rfloor - \frac{1}{2} \lfloor  x \rfloor\\
&= \theta \frac{x^2}{2} -x \{ \theta x \} + \frac{x}{2} + \frac{ \{
\theta x \}^2}{2 \theta} - \frac{ \{
\theta x \}}{2 \theta} +\\
& \quad + \theta \frac{x^2}{2} -\theta x \{ x \} + \theta \frac{x}{2} +\theta
\frac{ \{  x \}^2}{2} - \theta
\frac{ \{  x \}}{2} +\\  
& \quad - \bigl ( x - \{  x \} \bigr ) \bigl ( \theta x - \{ \theta x
\} \bigr ) - \frac{1}{2} \bigl ( \theta x - \{ \theta x
\} \bigr ) - \frac{1}{2} \bigl ( x - \{  x \} \bigr ).
\end{align*}

On constate que tous les termes o\`u $x$ appara\^{\i}t sans accolades
$\{ \, \}$ s'\'eliminent. Il reste la quantit\'e suivante :
$$
\frac{\theta}{2} \{ x \}^2 + \frac{1}{2 \theta} 
\{ \theta x \}^2 - \{ x \} \{ \theta x \}+ \frac{\theta -1}{2
\theta}\{ \theta x \} + \frac{1-\theta}{2}\{ x \},
$$
qui est bien \'egale \`a celle apparaissant dans l'\'enonc\'e.

Si $\theta$ est rationnel, $\theta=p/q$ avec $p \in \Nat^*$, $q \in
\Nat^*$, $(p,q)=1$, on a
$$
n/\theta \in \Int \quad \text{et} \quad 1 \ioe n \ioe \theta x \Ssi
\frac{n}{p} \in \Int \quad \text{et} \quad  1 \ioe \frac{n}{p} \ioe
\frac{x}{q},
$$
et
$$
m\theta \in \Int \quad \text{et} \quad 1 \ioe m \ioe  x \Ssi
\frac{m}{q} \in \Int \quad \text{et} \quad  1 \ioe \frac{m}{q} \ioe
\frac{x}{q}.
$$

Il faut donc ajouter $\lfloor x/q \rfloor$ au second membre de la
premi\`ere ligne du calcul ci-dessus. Mais en vertu de la proposition
\ref{t86}, on retranchera $\lfloor x/q \rfloor$ lors de l'\'etape
suivante, ce qui fait que le r\'esultat est inchang\'e.\fin

\begin{prop}\label{t110}
Soit $x$ et $\theta$ deux nombres r\'eels positifs. On a
\begin{multline*}
\sum_{1 \ioe m \ioe x}
\frac{B_1(m \theta)}{m} + \theta  \sum_{1 \ioe n \ioe \theta x}
\frac{B_1(n/\theta)}{n} = \frac{\theta}{2} \int_0^x  \{ t \}^2 t^{-2}
dt + \frac{1}{2} \int_0^{\theta x} \{ t \}^2 t^{-2}
dt - \int_0^x  \{ t \}\{ \theta t \} t^{-2}dt + \\
+ \frac{\theta - 1}{2}
\log \frac{1}{\theta} + \frac{\theta - 1}{2} \int_x^{\theta x}  \{ t
\} t^{-2} + \frac{1}{2\theta x} \bigl ( \{ \theta x \}-\theta\{ x \}
\bigr )^2 + \frac{\theta - 1}{2 \theta x} \bigl ( \{ \theta x \}-\theta\{ x \}
\bigr ).
\end{multline*}
\end{prop}
\dem

Posons
$$
S_{\theta}(x) := \sum_{1 \ioe m \ioe x} B_1(m \theta).
$$

La somme \`a \'evaluer peut s'\'ecrire au moyen d'int\'egrales de
Stieltjes :
\begin{align*}
\int_0^x t^{-1} dS_{\theta}(t) + \theta \int_0^{\theta x} t^{-1}
dS_{1/\theta}(t) &= \int_0^x u^{-1} d \bigl (S_{\theta}(u) +
S_{1/\theta}(\theta u) \bigr ) \\
&= \frac{S_{\theta}(x) + S_{1/\theta}(\theta x)}{x} + \int_0^x u^{-2}  \bigl (S_{\theta}(u) +
S_{1/\theta}(\theta u) \bigr ) du.
\end{align*}

D'apr\`es la proposition \ref{t88}, on a
$$
S_{\theta}(x) + S_{1/\theta}(\theta x) = \frac{1}{2\theta} \bigl ( \{ \theta x \}-\theta\{ x \}
\bigr )^2 + \frac{\theta - 1}{2 \theta} \bigl ( \{ \theta x \}-\theta\{ x \}
\bigr ).
$$

Or,
\begin{align*}
\int_0^x u^{-2}  \bigl (\{ \theta u \}-\theta\{ u \}
\bigr )^2 du &= \int_0^x u^{-2} \{ \theta u \}^2 du - 2 \theta
\int_0^x u^{-2} \{ u \} \{ \theta u \} du + \theta^2 \int_0^x u^{-2}
\{ u \}^2 du\\
&= \theta \int_0^{\theta x} u^{-2} \{ u \}^2 du + \theta^2 \int_0^x u^{-2}
\{ u \}^2 du - 2 \theta
\int_0^x u^{-2} \{ u \} \{ \theta u \} du \, ;
\end{align*}
avec la proposition \ref{t109}, cela fournit le r\'esultat annonc\'e. \fin

\newpage

\section{Sur la fonction $\Gamma$}

Nous rappelons dans ce chapitre\footnote{Voir par exemple {\sc N. Nielsen},
Handbuch der Theorie der Gammafunktion.} quelques propri\'et\'es de la fonction
$\Gamma$, de son logarithme, et de sa d\'eriv\'ee logarithmique
$$
\psi(z) := \frac{\Gamma'}{\Gamma}(z).
$$
\begin{prop}
La fonction $\psi$ est m\'eromorphe dans le plan complexe. Ses p\^oles
sont les nombres entiers n\'egatifs ou nuls $-n$, $n \in \Nat$, et on
a pour tout autre nombre complexe $z$
$$
\psi(z)= -\gamma - \frac{1}{z} + z \sum_{n=1}^{+\infty}
\frac{1}{n(n+z)}, 
$$
o\`u $\gamma$ est la constante d'Euler\footnote{On peut retenir cette expression sous la
forme $\sum_{n \soe 0} \left(\frac{1}{n}-\frac{1}{n+z} \right)$, o\`u
le symbole $\frac{1}{0}$ doit \^etre interpr\'et\'e comme valant $-\gamma$.}.
\end{prop}

\begin{prop}
Pour $z \in \Com$, $-z \not \in \Nat$, on a 
$$
\psi'(z)=\sum_{n \soe 0} \frac{1}{(n+z)^2} .
$$
\end{prop}

\begin{prop}\label{t12}
Pour $z \in \Com$ et $n$ entier positif tels que $-nz \not \in \Nat$, on a
$$
\log \Gamma(nz)= - \frac{n-1}{2} \log 2 \pi + \left (nz - \frac{1}{2}
\right) \log n + \sum_{k=0}^{n-1}   \log \Gamma \left(z +
\frac{k}{n} \right );
$$
en particulier,
$$
\sum_{k=1}^{n}   \log \Gamma \left(\frac{k}{n} \right )=\frac{n-1}{2}
\log 2 \pi -\frac{1}{2}\log n.
$$
\end{prop}
\begin{prop}
Pour $z \in \Com$ et $n$ entier positif tels que $-nz \not \in \Nat$, on a
$$
\psi(nz) = \log n + \frac{1}{n} \sum_{k=0}^{n-1}   \psi \left(z +
\frac{k}{n} \right );
$$
en particulier,
$$
\sum_{k=1}^{n}\psi \left(\frac{k}{n} \right )=-n(\log n+ \gamma).
$$
\end{prop}

\begin{prop}\label{t13}
On a	
$$
\int_0^1 \log \Gamma (x) dx = \frac{1}{2} \log 2 \pi;
$$
$$
\int_0^1 x\psi(x)dx = -\frac{1}{2} \log 2 \pi;
$$
$$
\int_u^1 x\psi'(x)dx = \log \Gamma (u)-u \psi(u) -\gamma, \quad u>0;
$$
$$
\int_0^1 x^2\psi'(x)dx =\log 2 \pi -\gamma.
$$
\end{prop}
{\bf D\'emonstration}

La premi\`ere int\'egrale est l'int\'egrale de Raabe\footnote{{\sc
N. Nielsen}, loc. cit. \S 34,
formule (17).}.

Pour la deuxi\`eme, on a
$$
\int_0^1 x\psi(x)dx = x\log \Gamma (x) dx \Bigr\rvert_0^1 - \int_0^1
\log \Gamma (x) dx = -\frac{1}{2} \log 2 \pi.
$$

Ensuite,
$$
\int_u^1 x\psi'(x)dx = x\psi(x)\Bigr\rvert_u^1 - \int_u^1 \psi(x)dx =
\log \Gamma (u)-u \psi(u) -\gamma.
$$

Enfin,
\begin{equation*}
\int_0^1 x^2\psi'(x)dx = x^2\psi(x)\Bigr\rvert_0^1 - \int_0^1
2x\psi(x)dx =
-\gamma +\log 2 \pi.\fine
\end{equation*}

\begin{prop}
$$
\psi(1+z)=\psi(z)+ \frac{1}{z}.
$$
\end{prop}

\begin{prop}
Pour tout nombre entier naturel $N$, on a
$$
\psi(N+1)= H_N-\gamma,
$$
o\`u 
$$
H_N := \sum_{k=1}^N \frac{1}{k}.
$$
\end{prop}

\begin{prop}
$$
\psi(1-z)=\psi(z)+ \pi \cot \pi z.
$$
\end{prop}

\begin{prop}\label{t76}
Le d\'eveloppement de Laurent de $(2 \pi)^{-s} \Gamma(s)$ en $s=-1$
commence par
$$
2\pi  \bigl (- (s+1)^{-1} + \gamma -1 + \log 2\pi + \dots \bigr )
$$
\end{prop}
{\bf D\'emonstration}

D'une part,
\begin{align*}
(2 \pi)^{-s} &= 2\pi e^{-(s+1) \log 2\pi} \\
&=  2\pi \bigl (1 -(s+1) \log 2 \pi + \dots \bigr ).
\end{align*}

D'autre part,
\begin{align*}
\Gamma(s) &= \frac{\Gamma(s+2)}{s(s+1)} \\
&= \frac{-1}{1-(s+1)} \cdot \frac{1}{s+1} \bigl ( \Gamma(1) + (s+1)
\Gamma'(1) + \dots \bigr ) \\
&= - \bigl ( (s+1)^{-1} +1 + \dots \bigr ) \bigl ( 1 -\gamma (s+1) +
\dots \bigr )\\
&= \frac{-1}{s+1} + \gamma -1 + \dots
\end{align*}

Le r\'esultat en d\'ecoule par multiplication.\fin

Nous poserons
$$
J(z) := J_{1,1}(z,0) = \int_0^{+\infty} \frac{\{t\} - \frac{1}{2}}{t+z} dt.
$$
On peut exprimer la fonction $J$ \`a l'aide de la fonction $\Gamma$.
\begin{prop}
Pour $z \in \Com \setminus ]-\infty,0]$, on a
$$
J(z) = - \log \Gamma(z) + (z - \frac{1}{2}) \log z -z + \frac{1}{2}
\log( 2\pi).
$$
En particulier, pour tout nombre entier positif $N$, on a
$$
J(N)= - \log N! + (N+\frac{1}{2}) \log N - N + \frac{1}{2} \log( 2\pi).
$$
\end{prop}

\begin{prop}
On a pour tout $z \in \Com \setminus ]-\infty,0]$,
$$
|J(z)| \ioe \frac{1 + \pi \sqrt{2}}{12 \dist (z,]-\infty,0])}.
$$
\end{prop}

Pour $x \in \Real$ et $z \in \Com \setminus ]- \infty, -x]$, on pose
$$
J_{1,2}(z,x) := \int_x^{+\infty} \frac{B_1(t)}{(t+z)^2} dt.
$$

\begin{prop}
La fonction $z \mapsto J_{1,2}(z,x)$ est analytique pour $z \in \Com
\setminus ]- \infty, -x]$. On a, pour $x + \Re z >0$,
$$
|J_{1,2}(z,x) | \ioe \frac{1}{2(x + \Re z )}.
$$
\end{prop}
 
\begin{prop}
Pour $z \in \Com \setminus ]-\infty,0]$, on a
$$
J'(z) = -J_{1,2}(z,0) = -\psi(z)  + \log z -
\frac{1}{2z}.
$$

En particulier, pour tout nombre entier positif $N$, on a
$$
J'(N) = \log N + \gamma +\frac{1}{2N}- H_N,
$$
o\`u $H_N : = 1 + \frac{1}{2} + \dots + \frac{1}{N}$.
\end{prop}

\newpage

\section{Sommes de Vassiounine}

Soit $q$ un nombre entier positif. Nous avons
$$
1 +X + \cdots + X^{q-1} = \frac{1-X^q}{1-X}.
$$

En d\'erivant, on obtient
\begin{align*}
1 + 2X + \cdots + (q-1) X^{q-2} &= \frac{-qX^{q-1}(1-X)
+1-X^q}{(1-X)^2} \\
&= \frac{1-qX^{q-1}+(q-1)X^q}{(1-X)^2}.
\end{align*}

Par cons\'equent,
$$
\sum_{n=0}^{q-1} nX^n = \frac{X}{(1-X)^2}(1-qX^{q-1}
+(q-1)X^q).
$$
\begin{prop}\label{t69}
Soit $z$ une racine $q$-\`eme de l'unit\'e. On a
$$
\sum_{n=0}^{q-1} nz^n=
\begin{cases}
  \frac{q(q-1)}{2}&   \text{si $z=1$;}\\
  \frac{q}{z-1}&   \text{si $z \not = 1$.}
\end{cases}
$$
\end{prop}
{\bf D\'emonstration}

Si $z=1$, la somme en question vaut
$$
1 + \cdots + (q-1) = \frac{q(q-1)}{2}.
$$
Si $z \not = 1$, la formule ci-dessus donne
\begin{align*}
\sum_{n=0}^{q-1} nz^n &= \frac{z}{(1-z)^2}(1-qz^{q-1}
+q-1)\\
&= \frac{qz(1-z^{q-1})}{(1-z)^2} \\
&= \frac{q(z-1)}{(1-z)^2} \\
&= \frac{q}{z-1}.\fine
\end{align*}

Si $p$ et $q$ sont deux nombres entiers premiers entre eux, $q$
\'etant positif, nous d\'efinissons la {\bf somme de Vassiounine}\footnote{V. I. Vassiounine, {\it Sur un syst\`eme biorthogonal
   reli\'e \`a l'hypoth\`ese de Riemann} (en russe), Alg. i
   An. {\bf7} (1995), 118-135 ; traduction anglaise dans St-Petersburg
Math. J. {\bf7} (1996), 405-419.}
$V(p,q)$ par la formule
\begin{align*}
V(p,q) :&= \sum_{k=1}^{q-1} \left \{ \frac{kp}{q} \right \} \cot
\frac{k \pi}{q} \\
&= \frac{1}{q} \sum_{k=1}^{q-1} ( kp \bmod q )\cot
\frac{k \pi}{q} \\
&= \sum_{k=1}^{q-1} \frac{k}{q}\cot
\frac{k \overline{p}\pi}{q},
\end{align*}
o\`u $p\overline{p} \equiv 1\pmod q$. Par convention, $V(p,1)=0$.

\begin{prop}\label{t75}
La somme de Vassiounine $V(p,q)$ est une fonction impaire et de
p\'eriode $q$ de la variable $p$ :
$$
V(p+q,q)=V(p,q); \quad V(-p,q) = -V(p,q).
$$
\end{prop}

\begin{prop}\label{t100}
On a
\begin{align*}
V(p,q) &= \sum_{k=1}^{q-1} B_1 \left ( \frac{kp}{q}
\right )\cot \frac{k \pi}{q} \\
&= 2 \sum_{1 \ioe k < q/2} B_1 \left ( \frac{kp}{q}
\right )\cot \frac{k \pi}{q}
\end{align*}
\end{prop}
{\bf D\'emonstration}

Comme $\cot (\pi -x) = - \cot x$, on a
\begin{align*}
V(p,q) &= \sum_{k=1}^{q-1} \left \{ \frac{kp}{q} \right \} \cot
\frac{k \pi}{q} \\
&= \sum_{k=1}^{q-1} \left \{ \frac{(q-k)p}{q} \right \} \cot
\frac{(q-k) \pi}{q} \\
&= - \sum_{k=1}^{q-1} \left \{- \frac{kp}{q} \right \} \cot
\frac{k \pi}{q} \\
&= \sum_{k=1}^{q-1} \frac{\{kp/q\}-\{-kp/q\}}{2} \cot
\frac{k \pi}{q} \\
&= \sum_{k=1}^{q-1} B_1 \left ( \frac{kp}{q}
\right )\cot \frac{k \pi}{q}.
\end{align*}

D'autre part, pour $1 \ioe k \ioe q-1$,
$$
 B_1 \left ( \frac{kp}{q} \right )\cot \frac{k \pi}{q} = B_1 \left
 ( \frac{(q-k)p}{q} \right )\cot \frac{(q-k) \pi}{q},
$$
donc la contribution \`a la derni\`ere somme de l'intervalle $q/2 < k
\ioe q-1$ est \'egale \`a celle de l'intervalle $1 \ioe k < q/2$
(celle de $k=q/2$ est nulle).\fin

\begin{prop}\label{t101}
Pour $p,q \in \Int$, $q \soe 2$, on a $V(p,q) \ll q \log q$.
\end{prop}
{\bf D\'emonstration}

Pour $0 < x \ioe \pi/2$, on a $\cot x = 1/x + O(1)$, d'o\`u
\begin{align*}
V(p,q) &=  \frac{2}{\pi} \sum_{1 \ioe k < q/2} \frac{B_1(kp/q)}{k} +
O(q) \\
& \ll q \log q.\fine
\end{align*}

On peut exprimer les sommes de Vassiounine au moyen de la d\'eriv\'ee
logarithmique de la fonction $\Gamma$.

\begin{prop}\label{t111}
On a
\begin{align*}
V(p,q) &= - \frac{2}{\pi} \sum_{k=1}^{q-1} B_1 \left ( \frac{kp}{q}
\right ) \psi \left ( \frac{k}{q} \right ) \\
&= - \frac{2}{\pi} \sum_{k=1}^{q-1} \left \{ \frac{kp}{q} \right \}
\psi \left ( \frac{k}{q} \right ) -\frac{q}{\pi} (\log q + \gamma ).
\end{align*}
\end{prop}
{\bf D\'emonstration}
\begin{align*}
V(p,q) &= \sum_{k=1}^{q-1} \left \{ \frac{kp}{q} \right \} \cot
\frac{k \pi}{q} \\
&= \frac{1}{\pi} \sum_{k=1}^{q-1} \left \{ \frac{kp}{q} \right \} \left ( \psi \left
( 1-\frac{k}{q} \right )- \psi \left ( \frac{k}{q} \right ) \right )\\
&= \frac{1}{\pi} \sum_{k=1}^{q-1} \left \{ \frac{kp}{q} \right \}  \psi \left
( \frac{q-k}{q} \right )-\frac{1}{\pi} \sum_{k=1}^{q-1} \left \{ \frac{kp}{q} \right \}  \psi \left
( \frac{k}{q} \right )\\
&= \frac{1}{\pi} \sum_{k=1}^{q-1} \left \{ \frac{(q-k)p}{q} \right \}  \psi \left
( \frac{k}{q} \right )-\frac{1}{\pi} \sum_{k=1}^{q-1} \left \{ \frac{kp}{q} \right \}  \psi \left
( \frac{k}{q} \right )\\
&= \frac{1}{\pi} \sum_{k=1}^{q-1} \left ( \left \{ \frac{-kp}{q}
\right \}- \left \{ \frac{kp}{q}
\right \} \right ) \psi \left
( \frac{k}{q} \right )\\
&= - \frac{2}{\pi} \sum_{k=1}^{q-1} B_1 \left ( \frac{kp}{q}
\right ) \psi \left ( \frac{k}{q} \right ) \\
&= - \frac{2}{\pi} \sum_{k=1}^{q-1} \left \{ \frac{kp}{q} \right \}
\psi \left ( \frac{k}{q} \right ) +\frac{1}{\pi} \sum_{k=1}^{q-1}\psi \left ( \frac{k}{q} \right )\\  
&= - \frac{2}{\pi} \sum_{k=1}^{q-1} \left \{ \frac{kp}{q} \right \}
\psi \left ( \frac{k}{q} \right ) -\frac{q}{\pi} (\log q + \gamma ).\fine
\end{align*}

On peut consid\'erer les sommes apparaissant dans la proposition
\ref{t111} en y rempla\c{c}ant $p$ et $q$ par des nombres entiers
positifs $a$ et $b$ non n\'ecessairement premiers entre eux.

\begin{prop}\label{t112}
Soit $a$ et $b$ des nombres entiers positifs, $d$ leur plus grand
diviseur commun, $p=a/d$ et $q=b/d$. On a :
$$
\sum_{m=1}^{b-1} B_1 \left ( \frac{ma}{b}
\right ) \psi \left ( \frac{m}{b} \right ) =-\frac{\pi d}{2} V(p,q),
$$
et
$$
\sum_{m=1}^{b-1}  \left \{ \frac{ma}{b}
\right ) \psi \left ( \frac{m}{b} \right \} =-\frac{\pi d}{2} V(p,q)-
\frac{b}{2} (\log b + \gamma) + \frac{d}{2} (\log d + \gamma).
$$
\end{prop}
\dem

Observons que l'on peut remplacer les $\sum_{m=1}^{b-1}$ par des
$\sum_{m=1}^b$.

Pour $1 \ioe m \ioe b =qd$, \'ecrivons
$$
m=qk+r; \quad 0 \ioe k \ioe d-1; \quad 1 \ioe r \ioe q.
$$

On a :
\begin{align*}
\sum_{m=1}^b B_1 \left ( \frac{ma}{b}
\right ) \psi \left ( \frac{m}{b} \right ) &= \sum_{k=0}^{d-1}
\sum_{r=1}^q  B_1 \left ( \frac{(qk+r)pd}{qd} \right ) \psi \left
( \frac{qk+r}{qd}  \right ) \\
&=  \sum_{k=0}^{d-1}
\sum_{r=1}^q  B_1 \left ( \frac{rp}{q} \right ) \psi \left
(\frac{k}{d} + \frac{r}{qd} \right )  \\
&= \sum_{r=1}^q  B_1 \left ( \frac{rp}{q} \right )\sum_{k=0}^{d-1}\psi
\left 
(\frac{k}{d} + \frac{r}{qd} \right )  \\
&= \sum_{r=1}^q  B_1 \left ( \frac{rp}{q} \right ) \left ( \psi \left (
\frac{r}{q} \right ) - \log d \right ) d \\
&= d \sum_{r=1}^q  B_1 \left ( \frac{rp}{q} \right ) \psi \left (
\frac{r}{q} \right ) \\
&= -\frac{\pi d}{2} V(p,q).
\end{align*}

D'autre part,
\begin{align*}
\sum_{m=1}^b  \left \{ \frac{ma}{b}
\right ) \psi \left ( \frac{m}{b} \right \} &=\sum_{m=1}^b B_1 \left ( \frac{ma}{b}
\right ) \psi \left ( \frac{m}{b} \right ) + \frac{1}{2} \sum_{m=1}^b
\psi \left ( \frac{m}{b} \right ) - \frac{1}{2} \sum_{k=1}^d \psi
\left ( \frac{k}{d} \right ) \\
&= -\frac{\pi d}{2} V(p,q) - \frac{b}{2} (\log b + \gamma) +
\frac{d}{2} (\log d + \gamma). \fine
\end{align*}

Nous \'evaluons maintenant d'autres sommes trigonom\'etriques au moyen
de sommes de Vassiounine.
\begin{prop}\label{t70}
Soit $p$ et $q$ deux nombres entiers premiers entre eux, $q$
\'etant positif. On a
$$
\sum_{1 \ioe k,l \ioe q} kl e^{i2 \pi klp/q} = \frac{q^2}{4} \bigl (3q
+1 -2i V(\overline{p},q) \bigr ).
$$
\end{prop}
{\bf D\'emonstration}

\begin{align*}
\sum_{1 \ioe k,l \ioe q} kl e^{i2 \pi klp/q} &= \sum_{1 \ioe k,l \ioe
q-1} + \sum_{\substack{1 \ioe l \ioe q-1 \\
             k=q}} + \sum_{\substack{1 \ioe k \ioe q-1 \\
             l=q}} +q^2 \\
&= \sum_{k=1}^{q-1} k \sum_{l=1}^{q-1} l e^{i2 \pi klp/q} + 2q
\sum_{j=1}^{q-1}j +q^2 \\
&= \sum_{k=1}^{q-1} k \frac{q}{e^{i2 \pi kp/q}-1} + 2q
\frac{q(q-1)}{2} + q^2 && \text{d'apr\`es la proposition \ref{t69}}\\
&= q \sum_{k=1}^{q-1} \frac{k}{e^{i2 \pi kp/q}-1} +q^3 \\
&= q \sum_{k=1}^{q-1} k \frac{e^{-i \pi kp/q}}{2i \sin \pi kp/q} +q^3
\\
&= q \sum_{k=1}^{q-1} k \left ( -\frac{i}{2} \cot \pi \frac{kp}{q} -
\frac{1}{2} \right ) + q^3 \\
&= -\frac{i}{2}q^2 \sum_{k=1}^{q-1}\frac{k}{q} \cot \pi \frac{kp}{q} -
\frac{1}{2} q \frac{q(q-1)}{2} + q^3 \\
&= \frac{3}{4} q^3 + \frac{1}{4} q^2 -\frac{i}{2}q^2
V(\overline{p},q),
\end{align*}
comme annonc\'e.\fin

\begin{prop}\label{t72}
Soit $p$ et $q$ deux nombres entiers premiers entre eux, $q$
\'etant positif. On a
$$
\sum_{1 \ioe k,l \ioe q} \left ( \frac{1}{2} - \frac{k}{q} \right )\left
( \frac{1}{2} - \frac{l}{q} \right )  e^{i2 \pi klp/q} =\frac{1}{4} - \frac{i}{2}
V(\overline{p},q).
$$
\end{prop}
{\bf D\'emonstration}

On a
\begin{multline*}
\sum_{1 \ioe k,l \ioe q} \left ( \frac{1}{2} - \frac{k}{q} \right )
\left( \frac{1}{2} - \frac{l}{q} \right )  e^{i2 \pi klp/q} =
\frac{1}{4}\sum_{1 \ioe k,l \ioe q}e^{i2 \pi klp/q} -
\frac{1}{2q}\sum_{1 \ioe k,l \ioe q}l e^{i2 \pi
klp/q}+\\
-\frac{1}{2q}\sum_{1 \ioe k,l \ioe q}k e^{i2 \pi klp/q} +
\frac{1}{q^2} \sum_{1 \ioe k,l \ioe q}kl e^{i2 \pi klp/q}.
\end{multline*}

La premi\`ere somme vaut
$$
\frac{1}{4}\sum_{\substack{1 \ioe l \ioe q\\
                           k=q}}1 + \frac{1}{4} \sum_{k=1}^{q-1}
                           \sum_{l=1}^q (e^{i2 \pi kp/q})^l =
                           \frac{q}{4}.
$$

La deuxi\`eme, comme la troisi\`eme, vaut
$$
-\frac{1}{2q} \left (\sum_{\substack{1 \ioe k \ioe q\\
                           l=q}}q + \sum_{l=1}^{q-1} l\sum_{k=1}^q
                           (e^{i2 \pi lp/q})^k\right ) = -
                           \frac{q}{2}.
$$

Par cons\'equent, la proposition \ref{t70} nous donne
\begin{align*}
\sum_{1 \ioe k,l \ioe q} \left ( \frac{1}{2} - \frac{k}{q} \right )
\left ( \frac{1}{2} - \frac{l}{q} \right )  e^{i2 \pi klp/q} &=- \frac{3}{4}
q + \frac{1}{q^2} \sum_{1 \ioe k,l \ioe q} kl e^{i2 \pi klp/q}\\
&= - \frac{3}{4}
q +\frac{3}{4}
q +\frac{1}{4}- \frac{i}{2}V(\overline{p},q)\\
&= \frac{1}{4}- \frac{i}{2}V(\overline{p},q).\fine
\end{align*}

\newpage

\section{Somme des inverses des termes d'une progression
arithm\'etique}

La source de ce paragraphe est l'article de {\sc D. H. Lehmer},
\textit{Euler constants for arithmetical progressions}, Acta
Arith. {\bf 27} (1975), 125--142.

Soient $q$ et $r$ des nombres entiers positifs.
\begin{prop}
On a pour $x>0$,
$$
\sum_{0 \ioe n \ioe x} \frac{1}{qn+r} = \frac{\log x}{q} -\frac{1}{q}\psi\left(\frac{r}{q}\right) + \frac{1}{q}\log
\left(1+ \frac{r}{qx} \right) + \frac{\frac{1}{2}- \{x\}}{qx+r}
+\frac{1}{q} J_{1,2}\left(\frac{r}{q},x \right).
$$
\end{prop}
{\bf D\'emonstration}

On a 
\begin{eqnarray*}
\sum_{0 \ioe n \ioe x} \frac{1}{qn+r} & = & \int_{0^-}^x \frac{d
\lfloor t \rfloor}{qt+r} \\
                                      & = & \int_0^x \frac{dt}{qt+r} +
\int_{0^-}^x \frac{d\left(\frac{1}{2}- \{t\}\right)}{qt+r}\\
                                      & = & \frac{1}{q}\log
\frac{qx+r}{r} + \frac{\frac{1}{2}- \{x\}}{qx+r} + \frac{1}{2r} +q
\int_0^x \frac{\frac{1}{2}- \{t\}}{(qt+r)^2} dt \\
                                      & = & \frac{\log x}{q} +
\frac{1}{q}\log \frac{q}{r} + \frac{1}{2r} + q\int_0^{+\infty}
\frac{\frac{1}{2}- \{t\}}{(qt+r)^2} dt + \frac{1}{q}\log
\left(1+ \frac{r}{qx} \right) + \frac{\frac{1}{2}- \{x\}}{qx+r}
-q\int_x^{+\infty} \frac{\frac{1}{2}- \{t\}}{(qt+r)^2} dt,
\end{eqnarray*}
d'o\`u le r\'esultat puisque
\begin{equation*}
\frac{1}{q}\log \frac{q}{r} + \frac{1}{2r} +
\frac{1}{q}J'\left(\frac{r}{q}\right) = -\frac{1}{q}  \psi
\left(\frac{r}{q} \right).\fine
\end{equation*}

\begin{prop}
On a pour $x>r$,
$$
\sum_{\substack{
         n \ioe x\\
         n \equiv r \bmod{q}}}
\frac{1}{n} = \frac{\log x}{q} + \gamma(r,q) + R(x,r,q),
$$
o\`u
$$ 
\gamma(r,q) := -\frac{1}{q} \left( \psi
\left(\frac{r}{q} \right) + \log q \right),
$$
et
$$
R(x,r,q):= \frac{\frac{1}{2}- \{\frac{x-r}{q}\}}{x}
+\frac{1}{q} J_{1,2}\left(\frac{r}{q},\frac{x-r}{q}\right).
$$
 
On a $R(x,r,q) \ioe 1/x$.
\end{prop}
{\bf D\'emonstration}

On a 
\begin{align*}
\sum_{\substack{
         n \ioe x\\
         n \equiv r \bmod{q}}}
\frac{1}{n}                     & =  \sum_{0 \ioe n \ioe (x-r)/q}
         \frac{1}{qn+r}\\
& = \frac{1}{q}\log \left(\frac{x-r}{q} \right) -\frac{1}{q}  \psi
\left(\frac{r}{q} \right)+ \frac{1}{q}\log
\left(1+ \frac{r}{x-r} \right) + R(x,r,q)\\
                                & = \frac{\log x}{q}+\gamma(r,q) + R(x,r,q).\fine
\end{align*}

\begin{prop}\label{t133}
Soit $g: \Nat^* \vers \Com$ une fonction de p\'eriode $q$. On a, pour
$x>q$,
$$
\sum_{n \ioe x} \frac{g(n)}{n} = S(g) \log x + \gamma(g) + R(x,g),
$$
o\`u
$$
S(g):= \frac{1}{q} \sum_{r=1}^q g(r),
$$
$$
\gamma(g) := \sum_{r=1}^q g(r) \gamma(r,q),
$$
et 
$$
R(x,g):= \sum_{r=1}^q g(r) R(x,r,q).
$$
On a $|R(x,g)| \ioe x^{-1}\sum_{r=1}^q |g(r)|$.
\end{prop}
{\bf D\'emonstration}

On a 
\begin{align*}
\sum_{n \ioe x} \frac{g(n)}{n} & = \sum_{r=1}^q g(r)\sum_{\substack{
         n \ioe x\\
         n \equiv r \bmod{q}}}
\frac{g(n)}{n} \\
                               & =  \sum_{r=1}^q g(r) \left
         ( \frac{\log x}{q} + \gamma(r,q) + R(x,r,q) \right) \\
                               & =  S(g) \log x + \gamma(g) + R(x,g).\fine
\end{align*}

\begin{prop}\label{t134}
Soit $g: \Int \vers \Com$ une fonction de p\'eriode $q$. La s\'erie
$\sum_{n \soe 1} g(n)/n$ converge si et seulement si $S(g)=0$ et on a
alors
$$
\sum_{n \soe 1} \frac{g(n)}{n} =\gamma(g)=\sum_{r=1}^q g(r)
\gamma(r,q)= -\frac{1}{q}\sum_{r=1}^q g(r) \psi\left(\frac{r}{q} \right).
$$
\end{prop} 
\dem

Cela r\'esulte de la proposition \ref{t133}.\fin

Voici une application de cette proposition.

\begin{prop}\label{t80}
Soit $p$ et $q$ deux nombres entiers premiers entre eux, $q$ \'etant
positif. La s\'erie
$$
\sum_{k \soe 1} \frac{B_1(kp/q)}{k}
$$
converge et a pour somme $\frac{\pi}{2q} V(p,q)$.
\end{prop}
{\bf D\'emonstration}

La fonction $g : \, k \mapsto B_1(kp/q)$ est p\'eriodique, de
p\'eriode $q$ et
\begin{align*}
S(g) &= \frac{1}{q}\sum_{r=1}^q B_1(rp/q) \\
&= \frac{1}{q}B_1(rp)\\
&=0.
\end{align*}
La s\'erie est donc convergente et a pour somme 
\begin{align*}
\gamma(g) &= -\frac{1}{q} \sum_{r=1}^q g(r) \psi (r/q) \\
&= -\frac{1}{q} \sum_{r=1}^q B_1(rp/q) \psi (r/q) \\
&= \frac{\pi}{2q} V(p,q).\fine
\end{align*}

\smallskip
\begin{prop}
Soit $g: \Int \vers \Com$ une fonction de p\'eriode $q$. On a
$$
\sum_{n \soe 1} \frac{g(n)}{n(n+1)} = g(0) + \frac{1}{q}\sum_{r=1}^q
\bigl (g(r-1)-g(r) \bigr) \psi\left(\frac{r}{q} \right).
$$
\end{prop}
\dem

Comme $g$ est born\'ee, on a par sommation d'Abel :
$$
\sum_{n \soe 1} \frac{g(n)}{n(n+1)} = g(0) +\sum_{n \soe 1} \frac{h(n)}{n},
$$
o\`u $h(n):= g(n)-g(n-1)$. Or $h$ est de p\'eriode $q$ comme $g$ et, d'une part,
$$
S(h):= \frac{1}{q} \sum_{r=1}^q h(r)=\frac{1}{q} \bigl ( g(q)-g(0) \bigr )=0,
$$
d'autre part,
$$
\gamma(h)=\sum_{r=1}^q h(r)
\gamma(r,q)= -\frac{1}{q}\sum_{r=1}^q h(r) \psi\left(\frac{r}{q} \right).
$$
Le r\'esultat d\'ecoule donc de la proposition \ref{t134}.\fin

\newpage

\section{La fonction d'Estermann}

\subsection{D\'efinitions et propri\'et\'es fondamentales}
Soit $h$ et $k$ deux nombres entiers premiers entre eux, $k$ \'etant
positif. La {\bf fonction d'Estermann\footnote{{\sc M. Jutila}, Lectures on a method in the theory of exponential sums, \S 1.1}} $E(s;h/k)$ est d\'efinie pour
$\Re s >1$ par la s\'erie de Dirichlet absolument convergente
$$
E(s;h/k) := \sum_{n=1}^{+\infty} \tau(n) e^{i2 \pi n h/k}n^{-s}.
$$
Observons que $E(s;h/k)=E(s;h'/k)$ si $h \equiv h' \pmod k$.

On a en particulier $E(s;0/1)= \zeta(s)^2$. Plus g\'en\'eralement, on
peut exprimer $E(s;h/k)$ \`a l'aide de la fonction $\zeta$ d'Hurwitz.

\begin{prop}\label{t71}
$$
E(s;h/k)=k^{-2s} \sum_{1 \ioe j,l \ioe k}e^{i2 \pi jl h/k}
\zeta(s,j/k)\zeta(s,l/k).
$$
\end{prop}
\begin{prop}
La fonction $E(s;h/k)$ se prolonge m\'eromorphiquement au plan
complexe, avec un seul p\^ole, double, en $s=1$. La partie polaire de
$E(s;h/k)$ en $s=1$ est
$$
k^{-1} (s-1)^{-2} + k^{-1}(2 \gamma - 2 \log k) (s-1)^{-1}.
$$
\end{prop}

\begin{prop}\label{t84}
La fonction $E(s;h/k)$ v\'erifie l'\'equation fonctionnelle suivante
$$
E(s;h/k)=2 (2 \pi)^{2s-2} \Gamma^2(1-s) k^{1-2s} \bigl
(E(1-s;\overline{h}/k) - \cos \pi s E(1-s;-\overline{h}/k) \bigr )
$$
o\`u $h\overline{h} \equiv 1 \pmod k$ ( si $k=1$, on pose
 $E(s;\overline{h}/k)=\zeta^2(s)$).
\end{prop}
\addtocounter{footnote}{-8}
\begin{prop}\label{t74}
On a\footnote{Cf. theorem $1$ dans {\sc M. Ishibashi}, {\it The value of the Estermann zeta
functions at $s=0$},  Acta Arith. {\bf 73} (1995), 357-361.}
$$
E(0;h/k) = \frac{1}{4} - \frac{i}{2} V(\overline{h},k),
$$
o\`u $h\overline{h} \equiv 1 \pmod k$ et o\`u $V(\overline{h},k)$
d\'esigne la somme de Vassiounine
$$
V(\overline{h},k) = \sum_{j=1}^{k-1} \frac{j}{k} \cot \frac{jh
\pi}{k}.
$$
\end{prop}
{\bf D\'emonstration}

D'apr\`es la proposition \ref{t71}, on a
\begin{align*}
E(0;h/k) &= \sum_{1 \ioe j,l \ioe k}e^{i2 \pi jl h/k}
\zeta(0,j/k)\zeta(0,l/k) \\
&= \sum_{1 \ioe j,l \ioe k}\left ( \frac{1}{2} - \frac{j}{k} \right )\left
( \frac{1}{2} - \frac{l}{k} \right )  e^{i2 \pi jlh/k},
\end{align*}
car $\zeta(0,a)=\frac{1}{2}-a$. La proposition \ref{t72} permet de conclure.\fin

Comme la fonction $\zeta$ d'Hurwitz, la fonction d'Estermann est \`a
croissance polyn\^omiale, uniform\'ement dans toute bande
verticale. Pour sa fonction de Lindel\"of
$$
\mu(\sigma) := \limsup_{|\tau| \vers +\infty} \, \log
|E(\sigma+i\tau;h/k)| \, / \log |\tau|,
$$
on a les estimations suivantes
\begin{align*}
\mu(\sigma) &= 0, && \sigma \soe 1 &&& (\text{s\'erie de Dirichlet})\\
\mu(\sigma) &= 1- 2 \sigma,  && \sigma \ioe 0 &&& (\text{\'equation
fonctionnelle et formule de Stirling})\\
\mu(\sigma) &\ioe 1-\sigma, && 0 <\sigma <1 &&& (\text{convexit\'e}).
\end{align*}

\subsection{Les fonctions $\Esin$ et $\Ecos$}

Nous d\'efinissons une fonction d'Estermann {\og en cosinus
\fg} et une {\og en sinus \fg}.
$$
\Ecos (s;h/k) := \sum_{n=1}^{+\infty} \frac{\tau(n)}{n^s} \cos 2 \pi n
\frac{h}{k} = \frac{1}{2} E(s;h/k) +\frac{1}{2} E(s;-h/k);
$$
$$
\Esin (s;h/k) := \sum_{n=1}^{+\infty} \frac{\tau(n)}{n^s} \sin 2 \pi n
\frac{h}{k} = \frac{1}{2i} E(s;h/k) -\frac{1}{2i} E(s;-h/k).
$$

Observons que, pour $k=1$, $\Esin =0$ et $\Ecos = \zeta^2$.

\begin{prop}
La fonction $\Esin$ est enti\`ere. De plus, l'abscisse de convergence
de sa s\'erie de Dirichlet est inf\'erieure o\`u \'egale \`a $1/2$.
\end{prop}
{\bf D\'emonstration}

D'une part, $E(s;h/k)$ et $E(s;-h/k)$ ont m\^eme partie polaire en
$s=1$.

D'autre part, le th\'eor\`eme de Schnee-Landau et la valeur $\mu(0)=1$
pour la fonction de Lindel\"of de la fonction d'Estermann fournissent
l'assertion sur l'abscisse de convergence.\fin

\begin{prop}
La fonction $\Esin$ v\'erifie l'\'equation
fonctionnelle
$$
\Esin (s;h/k)=2 (2 \pi)^{2s-2} \Gamma^2(1-s) k^{1-2s} (1 + \cos \pi s
) \Esin (1-s;\overline{h}/k),
$$
ou, en posant
$$
\widetilde{\Esin} (s;h/k) := \sin \frac{\pi s}{2} \Gamma(s) (2
\pi/k)^{-s} \Esin (s;h/k),
$$
$$
\widetilde{\Esin} (s;h/k) = \widetilde{\Esin} (1-s;\overline{h}/k).
$$ 
\end{prop}
{\bf D\'emonstration}

Posons
$$
\chi(s) := 2 (2 \pi)^{2s-2} \Gamma^2(1-s) k^{1-2s}.
$$

D'apr\`es la proposition \ref{t84}, on a
$$
E(s;h/k) = \chi(s) \bigl
(E(1-s;\overline{h}/k) - \cos \pi s E(1-s;-\overline{h}/k) \bigr ),
$$ 
et
$$
E(s;-h/k) = \chi(s) \bigl
(E(1-s;-\overline{h}/k) - \cos \pi s E(1-s;\overline{h}/k) \bigr ).
$$ 

Par cons\'equent,
$$
\Esin (s;h/k)= \chi(s) \bigl ( \Esin (1-s;\overline{h}/k) -\cos \pi s
\Esin(1-s;-\overline{h}/k) \bigr ).
$$
Mais, la fonction sinus \'etant impaire, on a $\Esin (s;-h/k)=-\Esin
(s;h/k)$, donc
$$
\Esin (s;h/k) = \chi(s)(1 + \cos \pi s ) \Esin (1-s;\overline{h}/k).
$$

Comme
\begin{align*}
(1 + \cos \pi s )\chi(s) &= 2 \cos ^2 \frac{\pi s}{2} \cdot 2 (2
\pi)^{2s-2} \Gamma^2(1-s) k^{1-2s}\\
&= \left ( 2\cos \frac{\pi s}{2}\sin  \frac{\pi s}{2} \Gamma(1-s)
\right ) \frac{2\cos \frac{\pi s}{2}}{\sin  \frac{\pi s}{2}}
\Gamma(1-s) (2
\pi)^{2s-2}  k^{1-2s}\\ 
&= \left (\frac{\pi}{\Gamma(s)} \right ) \frac{2\cos \frac{\pi s}{2}}{\sin  \frac{\pi s}{2}}
\Gamma(1-s) (2
\pi)^{2s-2}  k^{1-2s}\\ 
&= \frac{\cos \frac{\pi s}{2} \Gamma(1-s)}{\sin  \frac{\pi s}{2}
\Gamma(s)} (2\pi/k)^{2s-1} \\
&= \frac{\sin  \frac{\pi (1-s)}{2} \Gamma(1-s) (2\pi/k)^{s-1}}{\sin
\frac{\pi s}{2} \Gamma(s) (2\pi/k)^{-s}},
\end{align*}
on obtient la forme sym\'etrique annonc\'ee.\fin

\begin{prop}\label{t81}
On a
$$
\Esin (0;h/k)= \frac{1}{2} V(\overline{h},k),
$$
et
$$
\Esin (1;h/k) = \sum_{n=1}^{+\infty}\frac{\tau (n)}{n} \sin 2
\pi n \frac{h}{k} = -\frac{\pi^2}{2k} V(h,k).
$$
\end{prop}
{\bf D\'emonstration}

On a
\begin{align*}
\Esin (0;h/k) &= \frac{1}{2i} E(0;h/k) -\frac{1}{2i} E(0;-h/k)\\
&= \Im E(0;h/k) \\
&=-\frac{1}{2} V(\overline{h},k).
\end{align*}

D'autre part, l'\'equation fonctionnelle de la fonction $\Esin$ donne
\begin{align*}
-\frac{1}{2} V(\overline{h},k) &= \Esin (0;h/k)\\
&= 2 (2 \pi)^{-2} \Gamma^2(1) k \cdot 2 \, \Esin(1;\overline{h}/k) \\
&= k \pi^{-2} \Esin(1;\overline{h}/k),
\end{align*}
d'o\`u
\begin{equation*}
\Esin(1;\overline{h}/k)= - \frac{\pi^2}{2k} V(\overline{h},k).\fine
\end{equation*}

\begin{prop}
La fonction $\Ecos$ est m\'eromorphe dans $\Com$ avec un seul p\^ole,
double, en $s=1$, o\`u sa partie polaire est
$$
k^{-1}(s-1)^{-2} + k^{-1} (2 \gamma - 2 \log k) (s-1)^{-1}.
$$
\end{prop}

\begin{prop}
La fonction $\Ecos$ v\'erifie l'\'equation fonctionnelle
$$
\Ecos (s;h/k)=2 (2 \pi)^{2s-2} \Gamma^2(1-s) k^{1-2s} (1 - \cos \pi s
) \Ecos (1-s;\overline{h}/k),
$$
ou, en posant
$$
\widetilde{\Ecos} (s;h/k) := \cos \frac{\pi s}{2} \Gamma(s) (2
\pi/k)^{-s} \Ecos (s;h/k),
$$
$$
\widetilde{\Ecos} (s;h/k) = \widetilde{\Ecos} (1-s;\overline{h}/k).
$$ 
\end{prop}

\begin{prop}
$$
\Ecos(0;h/k) = \frac{1}{4}.
$$
\end{prop}
{\bf D\'emonstration}

On a
\begin{align*}
\Ecos(0;h/k) &= \frac{E(0;h/k)+E(0;-h/k)}{2} \\
&= \Re E(0;h/k) \\
&= \frac{1}{4}.\fine
\end{align*}

\subsection{Les fonctions $G_0$ et $G_1$}

Nous
utiliserons les fonctions
$$
G_0(s;h/k) := \cos \frac{\pi s}{2}\, \Ecos(s;h/k) - \sin \frac{\pi s}{2}
\, \Esin(s;h/k)
$$
et
$$
G_1(s;h/k) := (2\pi)^{-s} \Gamma(s) G_0(s+2;h/k).
$$
Si $k=1$, on a $G_0(s;h)= \cos \frac{\pi s}{2} \zeta(s)^2$ et $G_1(s;h)= -(2\pi)^{-s} \Gamma(s)\cos \frac{\pi s}{2} \zeta(s+2)^2$.

\begin{prop}\label{t73}
La fonction $G_0$ est m\'eromorphe dans $\Com$ avec un seul p\^ole,
simple, en $s=1$, o\`u son d\'eveloppement de Laurent commence par
$$
-\frac{\pi}{2k} \bigl ( (s-1)^{-1} + 2 \gamma - 2 \log k - \pi V(h,k) \bigr ).
$$
\end{prop}
{\bf D\'emonstration}

Au voisinage de $s=1$, on a
\begin{align*}
G_0(s;h/k) &= \cos \frac{\pi s}{2} \Ecos(s;h/k) - \sin \frac{\pi s}{2}
\Esin(s;h/k)\\
&= \left ( - \frac{\pi}{2} (s-1) + O \bigl ( (s-1)^3 \bigr ) \right
) \left ( \frac{1}{k} (s-1)^{-2} + \frac{2 \gamma - 2 \log k }{k}
(s-1)^{-1} + O(1) \right ) +\\
& \quad \Bigl ( -1 + O \bigl ( (s-1)^2 \bigr ) \Bigr ) \left ( -
\frac{\pi^2}{2k} V(h,k) +O ( s-1) \right )\\
&= -\frac{\pi}{2k} (s-1)^{-1} + \frac{\pi^2 V(h,k) -2 \pi (\gamma-\log
k)}{2k} + O(s-1).\fine
\end{align*}

\begin{prop}\label{t77}
La fonction $G_0$ v\'erifie l'\'equation fonctionnelle
$$
G_0(s;h/k)=2 (2 \pi)^{2s-2} \Gamma^2(1-s) k^{1-2s} \sin \pi s \,
G_0(1-s;\overline{h}/k).
$$
\end{prop}
{\bf D\'emonstration}

Posons
$$
\chi(s) := 2 (2 \pi)^{2s-2} \Gamma^2(1-s) k^{1-2s}.
$$

D'apr\`es les \'equations fonctionnelles de $\Esin$ et $\Ecos$, on a
\begin{align*}
G_0(s;h/k) &=\cos \frac{\pi s}{2} \Ecos(s;h/k) - \sin \frac{\pi s}{2}
\Esin(s;h/k)\\
&= \chi(s) \left (\cos \frac{\pi s}{2} (1 - \cos \pi s
) \Ecos (1-s;\overline{h}/k)- \sin \frac{\pi s}{2}(1 + \cos \pi s
) \Esin (1-s;\overline{h}/k) \right )\\
&= \chi(s) \left (\cos \frac{\pi s}{2} \cdot 2 \sin^2 \frac{\pi s}{2}
\Ecos (1-s;\overline{h}/k)- \sin \frac{\pi s}{2}\cdot 2 \cos^2
\frac{\pi s}{2} \Esin (1-s;\overline{h}/k)\right )\\
&= \chi(s)\sin \pi s \left (\sin \frac{\pi s}{2}\Ecos
(1-s;\overline{h}/k) - \cos \frac{\pi s}{2}\Esin
(1-s;\overline{h}/k)\right )\\
&= \chi(s)\sin \pi s \left (\cos \frac{\pi (1-s)}{2}\Ecos
(1-s;\overline{h}/k) -\sin \frac{\pi (1-s)}{2}\Esin
(1-s;\overline{h}/k)\right )\\
&= \chi(s) \sin \pi s \,
G_0(1-s;\overline{h}/k).\fine
\end{align*}

\begin{prop}
On a
$$
G_0(0;h/k) = \frac{1}{4}.
$$
\end{prop}
{\bf D\'emonstration}
\begin{equation*}
G_0(0;h/k) = \Ecos (0;h/k) = \frac{1}{4}.\fine
\end{equation*}

\begin{prop}
La fonction $G_1(s;h/k)$ est m\'eromorphe dans $\Com$ avec des p\^oles
en $0,-1,-2, \dots$

La partie polaire en $s=-1$ est
$$
\frac{\pi^2}{k} \left ( \frac{1}{(s+1)^2} + \frac{1 + \gamma - 2 \log
k - \log 2 \pi - \pi V(h,k)}{s+1} + \dots \right ).
$$

Le p\^ole en $s=-2$ est simple, avec r\'esidu $\pi^2/2$.
\end{prop}
{\bf D\'emonstration}

D'apr\`es les propositions \ref{t76} et \ref{t73}, on a au voisinage
de $s=-1$,
$$
(2 \pi)^{-s} \Gamma(s) = 2\pi \left ( \frac{-1}{s+1} + \gamma -1 +
\log 2 \pi + \dots \right ),
$$
et
$$
G_0(s+2) = -\frac{\pi}{2k}\left ( \frac{1}{s+1} + 2\gamma - 2 \log k - \pi V(h,k)+ \dots
\right ).
$$
On trouve alors par multiplication la  partie polaire de $G_1$ en
$s=-1$ :
$$ 
-\frac{\pi^2}{k}\left ( \frac{-1}{(s+1)^2} + \frac{2 \log k - 2 \gamma +
 \gamma -1 + \log 2 \pi + \pi V(h,k)}{s+1} + \dots \right ),
$$
d'o\`u le r\'esultat annonc\'e.

Au voisinage de $s=-2$, on a $\Gamma(s) \sim \frac{1}{2(s+2)}$ et $(2\pi)^2
G_0(0;h/k)=\pi^2$.\fin

\begin{prop}\label{t95}
La fonction $G_1(s;h/k)$ v\'erifie l'\'equation fonctionnelle suivante
$$
G_1(s;h/k) = k^{-2s-3} \frac{(s+2)(s+3)}{s(s+1)}
G_1(-s-3;\overline{h}/k).
$$\end{prop}
{\bf D\'emonstration}

D'apr\`es la proposition \ref{t77}, on a
\begin{align*}
G_1(s;h/k) &= (2 \pi)^{-s} \Gamma(s) G_0(s+2;h/k)\\
&= (2 \pi)^{-s} \Gamma(s)2 (2 \pi)^{2(s+2)-2} \Gamma^2(1-s-2) k^{1-2(s+2)} \sin \pi (s+2) \,
G_0(1-s-2;\overline{h}/k)\\
&= 2 (2 \pi)^{s+2}\Gamma(s)\Gamma^2(-1-s)\sin \pi s \, k^{-2s-3}
G_0(-1-s;\overline{h}/k).
\end{align*}

Or,
\begin{align*}
\Gamma(s)\Gamma^2(-1-s)\sin \pi s &=
\Gamma(s)\frac{\Gamma(1-s)}{s(s+1)}(s+2)(s+3)\Gamma(-3-s)\sin \pi s\\
&=\pi \frac{(s+2)(s+3)}{s(s+1)} \Gamma(-3-s),
\end{align*}
donc
\begin{align*}
G_1(s;h/k) &=2 (2 \pi)^{s+2} \cdot \pi \frac{(s+2)(s+3)}{s(s+1)}
\Gamma(-3-s) k^{-2s-3}
G_0(-1-s;\overline{h}/k)\\
&= k^{-2s-3}\frac{(s+2)(s+3)}{s(s+1)} (2
\pi)^{s+3}\Gamma(-3-s)G_0(-1-s;\overline{h}/k)\\
&= k^{-2s-3}\frac{(s+2)(s+3)}{s(s+1)}G_1(-s-3;\overline{h}/k).\fine
\end{align*}

\begin{prop}\label{t96}
Le polyn\^ome g\'en\'eralis\'e
$$
\Res \bigl ( G_1(s;h/k)t^{-s}, -2 \bigr ) + \Res \bigl (
G_1(s;h/k)t^{-s}, -1 \bigr )
$$
est
$$
\frac{\pi^2}{2} t^2 - \frac{\pi^2}{k}t \bigl ( \log t + \pi V(h,k) + 2
\log k + \log 2 \pi - \gamma -1 \bigr ).
$$
\end{prop}
{\bf D\'emonstration}

Au voisinage de $s=-1$, on a
$$
G_1(s;h/k)t^{-s} = \frac{\pi^2}{k} \left ( \frac{1}{(s+1)^2} + \frac{1 + \gamma - 2 \log
k - \log 2 \pi - \pi V(h,k)}{s+1} + O(1) \right ) \cdot t \Bigl ( 1
-(s+1) \log t + O \bigl ( (s+1)^2 \bigr ) \Bigr ),
$$
d'o\`u
\begin{equation*}
\Res \bigl ( G_1(s;h/k)t^{-s}, -1 \bigr ) = \frac{\pi^2}{k}t \bigl (-\log t + 1 + \gamma - 2 \log
k - \log 2 \pi - \pi V(h,k) \bigr ).\fine
\end{equation*}

\begin{prop}\label{t98}
La fonction $G_1(s;h/k)$ est \`a croissance polyn\^omiale,
uniform\'ement dans toute bande verticale.
\end{prop}
{\bf D\'emonstration}

C'est en effet le cas pour les fonctions $\Gamma(s) \cos \frac{\pi
s}{2}$, $\Gamma(s) \sin \frac{\pi
s}{2}$, $\Ecos$ et $\Esin$.\fin

\noindent{\bf Remarque}

Si l'on avait d\'efini $G_0$ par la formule
$$
\cos \frac{\pi s}{2}\, \Ecos(s;h/k) + \sin \frac{\pi s}{2}
\, \Esin(s;h/k),
$$
les r\'esultats ci-dessus seraient les m\^emes, sauf \`a changer
$V(h,k)$ en son oppos\'e.

\newpage

\section{Les fonctions $\fhi_n$}

\subsection{D\'efinition de $\fhi_n$ pour $n \soe 2$}

Soit $n$ un nombre entier sup\'erieur ou \'egal \`a $2$, et $x$
un nombre r\'eel positif. On pose
$$
\fhi_n(x) := \sum_{k \soe 1} \frac{B_n(kx)}{k^n}.
$$

\begin{prop}
La fonction $\fhi_n$ est p\'eriodique, de p\'eriode $1$. Elle est de
classe ${\cal C}^{n-2}$ et, pour $n \soe 3$, on a
$\fhi'_n=n\fhi_{n-1}$.
\end{prop}
\begin{prop}\label{t25}
Pour $n \soe 2$ et $x \in \Real$, on a 
$$
\fhi_n(x) \ioe M_n \zeta(n) \ioe 6 n!(2\pi)^{-n}.
$$
\end{prop}
\begin{prop}
Pour $n \soe 2$ et $x \in \Real$, on a $\fhi_n(-x)=(-1)^n\fhi_n(x)$.
\end{prop}
\begin{prop}
Pour $n \soe 2$, la s\'erie de Fourier de $\fhi_n$ est
$$
-n! \sum_{m \not = 0} \frac{\tau(|m|)}{(i 2 \pi m)^n} e^{i 2 \pi m x},
$$
o\`u $\tau(|m|)$ d\'esigne le nombre de diviseurs positifs de $|m|$.
\end{prop}
 
\subsection{La fonction $\fhi_1$}

Commen\c{c}ons par une interversion formelle de sommations :
\begin{align*}
\sum_{k \soe 1} \frac{B_1(kx)}{k} &=-\frac{1}{\pi}\sum_{k \soe 1}
\frac{1}{k} \sum_{l \soe 1} \frac{\sin 2 \pi l kx}{l}\\
&=-\frac{1}{\pi}\sum_{k \soe 1, \, l \soe 1} \frac{\sin 2 \pi l kx}{lk}\\
&=-\frac{1}{\pi}\sum_{m \soe 1} \frac{\tau (m)}{m} \sin 2 \pi mx.
\end{align*}

Les propositions \ref{t80} et \ref{t81} montrent que cette
manipulation est justifi\'ee quand $x$ est un nombre rationnel $p/q$,
o\`u $p$ et $q$ sont entiers, premiers entre eux, et $q$ positif : les
deux s\'eries convergent alors et ont la m\^eme somme $\frac{\pi}{2q}
V(p,q)$.

\begin{prop}
Les s\'eries
$$
\sum_{k \soe 1} \frac{B_1(kx)}{k} \quad \text{et} \quad
-\frac{1}{\pi}\sum_{m \soe 1} \frac{\tau (m)}{m} \sin 2 \pi mx
$$
convergent presque partout et dans $L^2(0,1)$ vers une m\^eme
fonction.
\end{prop}
{\bf D\'emonstration}

L'assertion sur $L^2(0,1)$ r\'esulte de l'\'egalit\'e
$$
\sum_{k = 1}^K \frac{B_1(kx)}{k} = -\frac{1}{\pi}\sum_{m \soe 1}
m^{-1}\bigl (\sum_{\substack{d|m\\
d \ioe K}}1\bigr ) \sin 2 \pi mx,
$$
et de la convergence de la s\'erie
$$
\sum_{m \soe 1} \frac{\tau(m)^2}{m^2}.
$$

L'assertion {\og presque partout \fg} est d\^ue \`a Chowla et
Walfisz\footnote{Cf. Hilfssatz $14$ dans {\sc S. Chowla, A. Walfisz}, {\it \"Uber eine
Riemannsche Identit\"at}, Acta Arith. {\bf 1} (1936), 87-112.}.
\begin{flushright} $\Box$ \end{flushright}
\smallskip

Nous noterons
$$
\fhi_1(x) := \sum_{k \soe 1} \frac{B_1(kx)}{k} 
$$ en tout point de convergence de cette s\'erie\footnote{L'ensemble
de ces points de convergence est pr\'ecis\'ement connu ; il co\"{\i}ncide
d'ailleurs avec l'ensemble des points de convergence de la s\'erie de
Fourier (\ref{t82}). Cf. {\sc J. R. Wilton}, {\it An approximate
functional equation with application to a problem of diophantine
approximation},
J. reine angew. Math. {\bf 169} (1933), 219-237, et {\sc R. de la
Bret\`eche et G. Tenenbaum}, {\it S\'eries trigonom\'etriques \`a
coefficients arithm\'etiques}, \`a para\^{\i}tre au Journal d'analyse math\'ematique.}. La fonction
$\fhi_1$ est d\'efinie presque partout, et appartient \`a
$L^2(0,1)$. Sa s\'erie de Fourier est
\begin{equation}\label{t82}
-\frac{1}{\pi}\sum_{m \soe 1} \frac{\tau (m)}{m} \sin 2 \pi mx.
\end{equation}

\subsection{Majoration du module de continuit\'e de $\fhi_2$}

\begin{prop}
On a
$$
\fhi_2(x)= \int_0^x 2 \fhi_1(t) dt + \frac{\pi^2}{36}.
$$
En particulier, $\fhi_2$ est presque partout d\'erivable et $\fhi_2'=2
\fhi_1$ presque partout.
\end{prop}

Afin de majorer le module de continuit\'e de $\fhi_2$, donnons une
majoration des sommes partielles de la s\'erie de Fourier de
$\fhi_1$\footnote{Cf. formule ($25_{VI}$) dans {\sc A. Walfisz}, {\it \"Uber einige
trigonometrische Summen}, Math. Z. {\bf 33} (1931), 564-601.}

\begin{prop}\label{t91}
Pour $x \in \Real$ et $K \soe 2$, on a
$$
\sum_{1 \ioe k \ioe K} \frac{\tau (k)}{k} \sin kx \ll \log K.
$$
\end{prop}
{\bf D\'emonstration}

On a
\begin{align*}
\sum_{1 \ioe k \ioe K} \frac{\tau (k)}{k} \sin kx &= \sum_{ab \ioe K}
\frac{\sin abx}{ab} \\
&= \sum_{1 \ioe a \ioe K} \frac{1}{a} \sum_{b \ioe K/a} \frac{\sin
bax}{b}.
\end{align*}

La somme int\'erieure est born\'ee (proposition \ref{t89}), donc le
r\'esultat d\'ecoule de l'estimation
\begin{equation*}
\sum_{1 \ioe a \ioe K} \frac{1}{a} \ll \log K.\fine
\end{equation*}

\begin{flushright} $\Box$ \end{flushright}
\smallskip

Si $f : \, \Real \vers \Com$ est uniform\'ement continue, son module
de continuit\'e est la fonction de $\delta$ d\'efinie par
$$
\omega(\delta;f) := \sup_{|x-y| \ioe \delta} |f(x)-f(y)|.
$$
On a $\omega(\delta;f) = o(1)$ quand $\delta$ tend vers $0$.

\begin{prop}\label{t90}
Pour $0 < \delta \ioe 1/2$, on a $\omega(\delta;\fhi_2) \ll \delta
\log 1/ \delta$.
\end{prop}
{\bf D\'emonstration}

On peut supposer $\delta$ suffisamment petit.

Si $|x-y| \ioe \delta$ et $K \soe 2$, on a
\begin{align*}
\fhi_2(x)-\fhi_2(y) &= 2 \int_x^y \fhi_1(t) dt \\
&= - \frac{2}{\pi} \int_x^y  \sum_{k \ioe K} \frac{\tau(k)}{k} \sin 2
\pi kt  \, dt - \frac{2}{\pi} \int_x^y  \sum_{k > K} \frac{\tau(k)}{k} \sin 2
\pi kt  \, dt.
\end{align*}

La premi\`ere int\'egrale est $\ll \delta \log K$, d'apr\`es la
proposition \ref{t91}. L'in\'egalit\'e de Schwarz montre que le
carr\'e de la deuxi\`eme int\'egrale est
\begin{align*}
& \ll \delta \int_0^1 \left ( \sum_{k > K} \frac{\tau(k)}{k} \sin 2
\pi kt  \, dt \right ) ^2 dt \\
& \ll \delta \sum_{k > K} \frac{\tau(k)^2}{k^2},
\end{align*}
d'apr\`es l'\'egalit\'e de Parseval.

Comme
$$
\sum_{k > K} \frac{\tau(k)^2}{k^2} \ll \frac{\log^3 K}{K},
$$
on obtient
$$
\fhi_2(x)-\fhi_2(y) \ll \delta \log K + \delta^{1/2} K^{-1/2}
\log^{3/2} K.
$$

En choisissant 
$$
K = \frac{1}{\delta} \log \frac{1}{\delta},
$$
on aboutit au r\'esultat annonc\'e.\fin

\subsection{Transform\'ee de Mellin de $\fhi_2(p/q \, +t) - \fhi_2(p/q)$}

Soit $p$ et $q$ deux nombres entiers premiers entre eux, $q$ \'etant
positif. La fonction
$$
t \mapsto \Delta_{p,q}(t) := \fhi_2 \left (\frac{p}{q} +t \right ) - \fhi_2 \left (
\frac{p}{q} \right )
$$
appartient \`a $\Dbv (-1,0)$ en tant que fonction continue, born\'ee,
et $O(t \log 1/t)$ au voisinage de $0$ (proposition \ref{t90}). Si
$q=1$, $\Delta_{p,q}(t)=\fhi_2 (t) - \fhi_2 (0)$. Dans ce cas, la
fonction $\Delta_{\overline{p},q}(t)$, qui joue un r\^ole ci-dessous,
vaut par convention la m\^eme chose. Notons aussi que notre
d\'efinition de $\Delta_{p,q}(t)$ vaut pour tout nombre r\'eel $t$,
mais que, dans ce paragraphe, nous ne consid\'erons que les valeurs
positives de $t$.

\begin{prop}\label{t92}
Pour $-1 < \Re s < 0$, on a
$$
M \Delta_{p,q} (s) = -\frac{1}{\pi^2} G_1(s;p/q).
$$
\end{prop}
{\bf D\'emonstration}

On a
$$
\fhi_2(t) = \frac{1}{\pi^2} \sum_{m \soe 1} \frac{\tau(m)}{m^2} \cos 2
\pi mt,
$$
donc
\begin{align*}
\Delta_{p,q}(t) &= \frac{1}{\pi^2} \sum_{m \soe 1} \frac{\tau(m)}{m^2} \cos 2
\pi m \left ( \frac{p}{q} +t \right ) - \frac{1}{\pi^2} \sum_{m \soe 1} \frac{\tau(m)}{m^2} \cos 2
\pi m \frac{p}{q}\\
&= \frac{1}{\pi^2} \sum_{m \soe 1} \frac{\tau(m)}{m^2} \cos 2
\pi m \frac{p}{q} \cdot (\cos 2
\pi mt -1) - \frac{1}{\pi^2} \sum_{m \soe 1} \frac{\tau(m)}{m^2} \sin 2
\pi m \frac{p}{q} \cdot\sin 2
\pi mt.
\end{align*}

Comme $t \mapsto \cos t -1$ et $t \mapsto \sin t$ appartiennent \`a
$\Dbv (-1,0)$, et comme
$$
\sum_{m \soe 1} \frac{\tau(m)}{m^2} (2 \pi m)^{-\sigma} < + \infty
$$
pour tout $\sigma \in ]-1,0[$, la proposition \ref{t93} confirme que
$\Delta_{p,q} \in \Dbv (-1,0)$ et prouve que, pour $-1 < \Re s < 0$,
\begin{align*}
M\Delta_{p,q} (s) &= \frac{1}{\pi^2} \left ( \cos \frac{\pi s}{2}
\Gamma(s)  \sum_{m \soe 1} \frac{\tau(m)}{m^2} \cos 2
\pi m \frac{p}{q} \cdot (2 \pi m)^{-s} - \sin \frac{\pi s}{2}
\Gamma(s)  \sum_{m \soe 1} \frac{\tau(m)}{m^2} \sin  2
\pi m \frac{p}{q} \cdot (2 \pi m)^{-s} \right )\\
&= \frac{1}{\pi^2}\Gamma(s)(2 \pi )^{-s} \bigl ( \cos \frac{\pi s}{2}
\Ecos (s+2;p/q) - \sin \frac{\pi s}{2} \Esin(s+2;p/q) \bigr )\\
&=-\frac{1}{\pi^2}\Gamma(s)(2 \pi )^{-s} G_0(s+2;p/q)\\
&=-\frac{1}{\pi^2} G_1(s;p/q).\fine
\end{align*}

\begin{prop}\label{t99}
Soit $p$ et $q$ deux nombres entiers premiers entre eux, $q \soe 2$. La fonction
$$
f \, : \, t \mapsto  (qt)^3  \Delta_{\overline{p},q}(1/q^2t) - 2q^3
\int_0^t u(3u-t)\Delta_{\overline{p},q}(1/q^2u) du 
$$
appartient \`a $\Dbv (-3,-2)$ et v\'erifie
$$
Mf(s) = -\frac{1}{\pi^2} G_1(s;p/q), \quad -3 < \sigma < -2.
$$
\end{prop}
{\bf D\'emonstration}

Nous allons appliquer syst\'ematiquement les r\`egles de calcul de
transform\'ees de Mellin vues au chapitre \ref{t94}.

$\bullet$ La fonction
$$
f_1  \, : \, t \mapsto \Delta_{\overline{p},q}(t)
$$
appartient \`a $\Dbv (-1,0)$ et v\'erifie
$$
Mf_1 (s) = -\frac{1}{\pi^2} G_1(s;\overline{p}/q), \quad -1 < \sigma < 0.
$$
\smallskip

$\bullet$ La fonction
$$
f_2  \, : \, t \mapsto \Delta_{\overline{p},q}(1/t)
$$
appartient \`a $\Dbv (0,1)$ et v\'erifie
\begin{align*}
Mf_2 (s) &= Mf_1 (-s) \\
&= -\frac{1}{\pi^2} G_1(-s;\overline{p}/q), \quad 0 < \sigma < 1.
\end{align*}
\smallskip

$\bullet$ La fonction
$$
f_3  \, : \, t \mapsto t^3  \Delta_{\overline{p},q}(1/t)
$$
appartient \`a $\Dbv (-3,-2)$ et v\'erifie
\begin{align*}
Mf_3 (s) &= Mf_2 (s+3)\\
&= -\frac{1}{\pi^2} G_1(-s-3;\overline{p}/q), \quad -3 < \sigma < -2.
\end{align*}
\smallskip

$\bullet$ La fonction
$$
f_4  \, : \, t \mapsto q^{-3} f_3(q^2t) = q^3 t^3  \Delta_{\overline{p},q}(1/q^2t)
$$
appartient \`a $\Dbv (-3,-2)$ et v\'erifie
\begin{align*}
Mf_4 (s) &= q^{-2s-3} Mf_3 (s)\\
&= -\frac{1}{\pi^2} q^{-2s-3}G_1(-s-3;\overline{p}/q), \quad -3 < \sigma < -2.
\end{align*}
\smallskip

$\bullet$ La fonction
\begin{align*}
f_ 5 \, : \, t & \mapsto f_4(t) -6 \int_0^t f_4(u) u^{-1} du + 2 t
\int_0^t f_4(u) u^{-2} du \\
&= q^3 \left (t^3 \Delta_{\overline{p},q}(1/q^2t) - 6 \int_0^t u^2
\Delta_{\overline{p},q}(1/q^2u) du + 2 t \int_0^t u
\Delta_{\overline{p},q}(1/q^2u) du \right )
\end{align*}
appartient \`a $\Dbv (-3,-2)$ et v\'erifie
\begin{align*}
Mf_5 (s) &= \left ( 1 + \frac{6}{s} - \frac{2}{s+1} \right ) Mf_4(s)\\
&= -\frac{1}{\pi^2} q^{-2s-3}
\frac{(s+2)(s+3)}{s(s+1)}G_1(-s-3;\overline{p}/q) \\
&=-\frac{1}{\pi^2} G_1(s;p/q), \quad -3 < \sigma <
-2,
\end{align*}
d'apr\`es la proposition \ref{t95}.\fin

\subsection{Comportement de $\fhi_2$ au voisinage d'un nombre rationnel}

Les propositions \ref{t97}, \ref{t96}, \ref{t98}, \ref{t92}, \ref{t99}
montrent que $\Delta_{p,q}$ v\'erifie l'\'equation fonctionnelle
suivante.

\begin{prop}
Pour $t > 0$ et $p$, $q$, deux nombres entiers premiers entre eux, $q
>0$, on a
$$
\Delta_{p,q} (t) = \frac{t \log t }{q}  + \frac{t}{q} \bigl ( \pi
V(p,q) + 2 \log q + \log 2 \pi - \gamma -1 \bigr )  - \frac{t^2}{2} + (qt)^3  \Delta_{\overline{p},q}(1/q^2t) - 2q^3
\int_0^t u(3u-t)\Delta_{\overline{p},q}(1/q^2u) du.
$$
\end{prop}

Pour la fonction $\Delta_{p,q} (-t)$, $t > 0$, des calculs analogues
m\`enent \`a l'\'equation fonctionnelle
$$
\Delta_{p,q} (-t) = \frac{t \log t }{q}  + \frac{t}{q} \bigl (- \pi
V(p,q) + 2 \log q + \log 2 \pi - \gamma -1 \bigr )  - \frac{t^2}{2} + (qt)^3  \Delta_{\overline{p},q}(-1/q^2t) - 2q^3
\int_0^t u(3u-t)\Delta_{\overline{p},q}(-1/q^2u) du.
$$

On en d\'eduit notamment le comportement asymptotique de $\fhi_2$ au
voisinage d'un nombre rationnel quelconque.
\begin{prop}\label{t102}
On a, uniform\'ement pour $t \in \Real$, $p \in \Int$, $q \in \Nat^*$,
$(p,q)=1$,
$$
\fhi_2 \left ( \frac{p}{q} + t \right ) - \fhi_2 \left (
\frac{p}{q}\right ) = \frac{1}{q} |t| \log |t| + \frac{t}{q}  \pi
V(p,q) + \frac{|t|}{q}\bigl (2 \log q + \log 2 \pi - \gamma -1 \bigr ) t - \frac{t^2}{2} +
O \bigl ( (qt)^3 \bigr ).
$$
\end{prop}

Observons que le terme 
$$
\frac{t}{q}  \pi
V(p,q) = 2 \fhi_1 \left ( \frac{p}{q} \right ) t
$$
correspond \`a la
d\'erivation terme \`a terme de la s\'erie de Fourier de $\fhi_2$.
\newpage

\section{L'espace de Hilbert $\Hil := L^2(0,+\infty;t^{-2}dt)$}

\subsection{Dilatations dans $\Hil$}

Consid\'erons l'espace de Hilbert suivant :
$$
\Hil :=  L^2(0,+\infty;t^{-2}dt).
$$

Rappelons\footnote{Cf. \S \ref{t113}.} que la transformation de
Mellin-Plancherel d\'efinit un op\'erateur unitaire entre
$\Hil$ et $L^2(-1/2 + \Real i; d\tau/2 \pi)$.

Pour tout nombre r\'eel positif $\lambda$, on d\'efinit la dilatation
$K_{\lambda} : \Hil \vers \Hil$ par
$$
K_{\lambda}f(x) = \lambda^{-1/2} f(\lambda x), \quad 0 < x < +\infty.
$$

\begin{prop}\label{t114}
L'application 
\begin{align*}
]0,+\infty[ \times \Hil & \vers \Hil \\
(\lambda,f) & \mapsto K_{\lambda}f
\end{align*}
est continue.
\end{prop}

\begin{prop}\label{t115}
Les $K_{\lambda}$ forment un groupe d'op\'erateurs unitaires de $\Hil$
:
$$
K_{\lambda}K_{\mu}=K_{\lambda \mu}; \quad
K_{1/\lambda}=K_{\lambda}^{-1}=K_{\lambda}^*.
$$
\end{prop}

\begin{prop}\label{t116}
Pour $f \in \Hil$ et $\lambda >0$, on a :
$$
MK_{\lambda}f(s) = \lambda^{-1/2-s} Mf(s).
$$
\end{prop}

\subsection{Fonctions d'autocorr\'elation}\label{t120}

Soit $\fhi$ un \'el\'ement de $\Hil$. Pour $\lambda >0$, on pose
\begin{align*}
A_{\fhi}(\lambda) & := \int_0^{+\infty} \fhi(t) \overline{\fhi(\lambda
t)} t^{-2}dt \\
&= \sqrt{\lambda} \langle \fhi,K_{\lambda}\fhi \rangle.
\end{align*}

\begin{prop}\label{t117}
L'application 
\begin{align*}
]0,+\infty[ \times \Hil & \vers \Com \\
(\lambda,\fhi) & \mapsto A_{\fhi}(\lambda)
\end{align*}
est continue.
\end{prop}

Observons que $|A_{\fhi}(\lambda)| \ioe \| \fhi \|^2 \sqrt{\lambda}$,
donc $A_{\fhi}$ se prolonge par continuit\'e en $\lambda=0$ en posant $A_{\fhi}(0)=0$.

\begin{prop}\label{t119} 
Pour $\lambda > 0$ et $\fhi \in \Hil$, 
$$
A_{\fhi}(\lambda) = \lambda\overline{A_{\fhi}(1/\lambda)}.
$$
\end{prop}
\dem

Cela r\'esulte du changement de variable $u=\lambda t$.\fin

\begin{prop}\label{t118}
Si $\fhi \in \Dbv (a,b)$, o\`u $a<-1/2<b$, alors $A_{\fhi} \in \Dbv
\bigl ( \max(a,-b-1), \min(b,-a-1) \bigr )$, et
$$
MA_{\fhi} (s) = M\fhi (-s-1)(M\fhi)^*(s).
$$
\end{prop}
\dem

Posons pour $t>0$,
$$
f(t) = t^{-1} \fhi(t) \quad \text{et} \quad g(t) = \overline{\fhi
(1/t)},
$$
de sorte que $A_{\fhi}(\lambda)=f*g(1/\lambda)$.

D'autre part,
\begin{align*}
f \in \Dbv(a+1,b+1) \quad &\text{et} \quad Mf(s)=  M\fhi (s-1);\\
g \in \Dbv(-b,-a) \quad &\text{et} \quad Mg(s) =(M\fhi)^*(-s).
\end{align*}

Par cons\'equent,
$$
f*g \in \Dbv \bigl ( \max (a+1,-b), \min (b+1,-a) \bigr )  \quad
\text{et} \quad Mf*g(s) = M\fhi (s-1)(M\fhi)^*(-s),
$$
d'o\`u le r\'esultat.\fin

\newpage

\section{L'autocorr\'elation multiplicative \\ de la fonction 
\og partie
fractionnaire \fg }

\subsection{D\'efinition et premi\`eres propri\'et\'es de la fonction
$A(\lambda)$}

La fonction d'autocorr\'elation multiplicative de la fonction {\og
partie fractionnaire \fg} est d\'efinie pour  $\lambda \soe 0$ par l'int\'egrale
$$
A(\lambda) := \int_0^{+\infty} \{t\}\{\lambda t\} \frac{dt}{t^2}.
$$

Observons pour commencer que $A(0)=0$ et $A(\lambda)>0$ pour tout
$\lambda>0$.

La fonction {\og
partie fractionnaire \fg} appartient \`a l'espace de
Hilbert $\Hil$, et $A$ n'est autre que sa fonction
d'autocorr\'elation, au sens du \S \ref{t120}. On a donc le r\'esultat
suivant.

\begin{prop}\label{t121}
$A(\lambda)$ est une fonction continue de $\lambda$ pour $\lambda \soe
0$, et
$$
A(\lambda) = \lambda A(1/\lambda), \quad \lambda >0.
$$
\end{prop}

D'autre part, la fonction {\og
partie fractionnaire \fg} appartient \`a $\Dbv(-1,0)$ et\footnote{{\sc
E. C. Titchmarsh}, The theory of the Riemann zeta function, formule (2.1.5).}
$$
\int_0^{+\infty} \{ t \} t^{s-1} dt = \frac{\zeta(-s)}{s}, \quad -1 <
\Re s < 0.
$$

La proposition \ref{t118} m\`ene donc au r\'esultat suivant.
\begin{prop}\label{t122}
La fonction $A$ appartient \`a $\Dbv(-1,0)$ et
$$
MA(s) = -\frac{\zeta(-s)\zeta(s+1)}{s(s+1)}.
$$\end{prop}

\subsection{Premi\`ere relation entre $A$ et $\fhi_1$}

\begin{prop}\label{t123}
$A(1) = \log 2 \pi - \gamma$.
\end{prop}
\dem

\begin{align*}
A(1) &= \int_0^{+\infty} \{ t \}^2 t^{-2} dt \\
&= \int_0^1 t^2 \sum_{n \soe 0} (n+t)^{-2} dt \\
&= \int_0^1 t^2 \psi' (t) dt \\
&= \log 2 \pi - \gamma. \fine
\end{align*}

\begin{prop}\label{t124}
Soit $\lambda$ un nombre r\'eel positif tel que la s\'erie $\fhi_1
(\lambda)=\sum_{k \soe 1} B_1(k \lambda)/k$ converge. Alors il en est
de m\^eme pour $\fhi_1 (1/\lambda)$ et
$$
A(\lambda) = \frac{1-\lambda}{2} \log \lambda + \frac{\lambda +1}{2}
(\log 2 \pi - \gamma) - \fhi_1 (\lambda) -\lambda\fhi_1 (1/\lambda).
$$
\end{prop}
\dem

En faisant tendre $x$ vers l'infini dans la proposition \ref{t110}, on
obtient :
$$
\fhi_1 (\lambda) +\lambda\fhi_1 (1/\lambda) = \frac{\lambda}{2} A(1) +
\frac{1}{2} A(1) - A(\lambda) + \frac{\lambda -1}{2} \log
\frac{1}{\lambda},
$$
d'o\`u le r\'esultat. \fin

En particulier, la fonction $\lambda \mapsto \fhi_1 (\lambda)
+\lambda\fhi_1 (1/\lambda)$ se prolonge en une fonction continue sur
$]0,+\infty[$. En choisissant $\lambda$ rationnel, on obtient le
r\'esultat suivant.
\begin{prop}\label{t125}
Soit $p$ et $q$ deux nombres entiers positifs premiers entre eux et
$\lambda=p/q$. On a
$$
A(\lambda)= \frac{1-\lambda}{2} \log \lambda + \frac{\lambda
+1}{2}(\log 2 \pi - \gamma) - \frac{\pi}{2q} \bigl ( V(p,q) + V(q,p)
\bigr ).
$$
\end{prop}

\subsection{Repr\'esentation graphique de la fonction $A(\lambda)$}

En calculant $A(\lambda)$ gr\^ace \`a la proposition \ref{t125} quand
$\lambda$ d\'ecrit la suite de Farey d'ordre 287, on obtient la
repr\'esentation graphique suivante.

$$
\epsfig{file=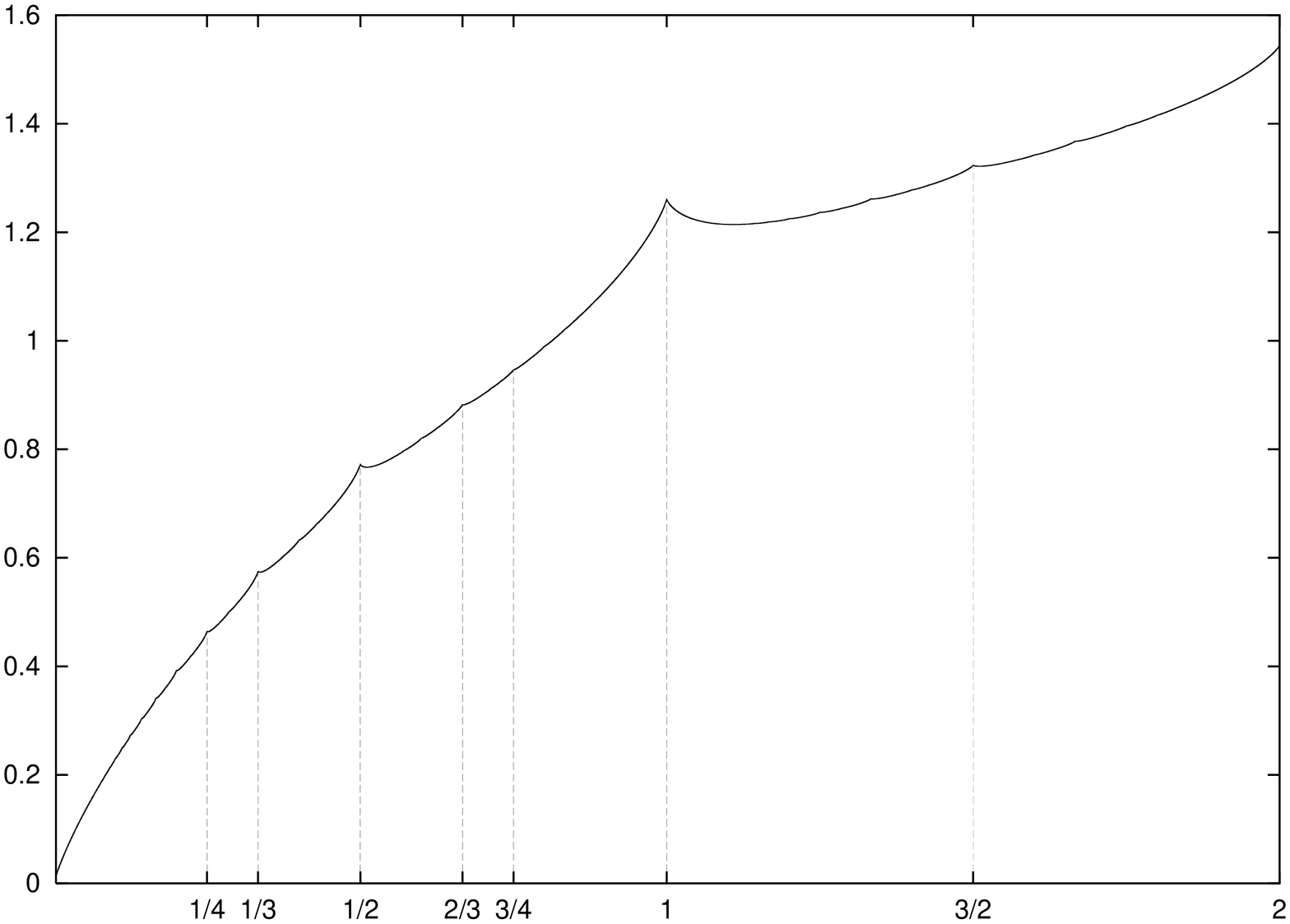,height=14.5cm,width=10cm,angle=270}
$$
\medskip

\subsection{Deuxi\`eme relation entre $A$ et $\fhi_1$}

\begin{prop}\label{t106}
On a pour tout $\lambda>0$,
\begin{align*}
A(\lambda) &= \frac{1}{2} \log \lambda + \frac{1-\gamma + \log 2\pi
}{2}-\lambda \int_{\lambda}^{+\infty}\fhi_1(t) \frac{dt}{t^2}\\
&= \frac{1}{2} \log \lambda + \frac{1-\gamma + \log 2\pi
}{2} + \frac{\fhi_2(\lambda)}{2 \lambda} - \lambda
\int_{\lambda}^{+\infty}\fhi_2(t) \frac{dt}{t^3}.
\end{align*}
\end{prop}
\dem

D'abord, une int\'egration par parties donne bien
$$
\int_{\lambda}^{+\infty}\fhi_1(t) \frac{dt}{t^2}=-\frac{\fhi_2(\lambda)}{2 \lambda^2} +\int_{\lambda}^{+\infty}\fhi_2(t) \frac{dt}{t^3}.
$$

D'autre part, pour tout $c \in \Com$ et tout $\lambda >0$, on a
$$
-\frac{\fhi_2(\lambda)}{2 \lambda^2}
 +\int_{\lambda}^{+\infty}\fhi_2(t) \frac{dt}{t^3} =
 -\frac{\fhi_2(\lambda)-c}{2 \lambda^2} +\int_{\lambda}^{+\infty}\bigl
 ( \fhi_2(t)-c\bigr ) \frac{dt}{t^3}.
$$

Choisissons $c= \fhi_2 (0)$, de sorte que $\fhi_2 (t) - c =
\Delta_{0,1}(t)$. D'apr\`es la proposition \ref{t92}, on sait que
$\Delta_{0,1} \in \Dbv(-1,0)$ et que
$$
M\Delta_{0,1} (s) = -\frac{1}{\pi^2} G_1(s;0).
$$

Par cons\'equent, la fonction
$$
f: \, t \mapsto \frac{\Delta_{0,1}(t)}{2t} -t \int_t^{+\infty}
\Delta_{0,1}(u) u^{-3} du
$$
appartient \`a $\Dbv(0,1)$ et v\'erifie
\begin{align*}
Mf(s) &= \frac{1}{2} M\Delta_{0,1} (s-1)- \frac{M\Delta_{0,1} (s-1)}{s+1}\\
&= -\frac{1}{2\pi^2} \frac{s-1}{s+1}G_1(s-1;0).
\end{align*}

V\'erifions que cette derni\`ere fonction n'est autre que
$$
-\frac{\zeta(-s)\zeta(s+1)}{s(s+1)}.
$$

On a 
\begin{align*}
-\frac{1}{2\pi^2} \frac{s-1}{s+1}G_1(s-1;0) &= -\frac{1}{s+1}
 \frac{s-1}{2\pi^2} \left ( -(2\pi)^{1-s} \Gamma(s-1)\cos \frac{\pi
 (s-1)}{2} \zeta(s+1)^2 \right )\\
&= - \frac{\zeta(s+1)}{s(s+1)}\cdot \Gamma(s+1)2 (2\pi)^{-1-s} \sin \frac{\pi
 (-s)}{2} \zeta (1+s)\\
&= -\frac{\zeta(-s)\zeta(s+1)}{s(s+1)},
\end{align*}
d'apr\`es l'\'equation fonctionnelle de la fonction $\zeta$ sous la
forme
$$
\zeta (s) = 2 (2\pi)^{s-1}\Gamma(1-s)\sin \frac{\pi
 s}{2} \zeta (1-s).
$$

Il r\'esulte alors des propositions \ref{t97} et \ref{t122} que
$$
f(\lambda) = A(\lambda) + \Res \left
(-\frac{\zeta(-s)\zeta(s+1)}{s(s+1)} \lambda^{-s}, 0 \right ),
$$
d'abord pour presque tout $\lambda >0$, mais en fait pour tout
$\lambda >0$ par continuit\'e de $f$ et $A$.

Au voisinage de $0$, on a les d\'eveloppements suivants :

\begin{align*}
\zeta(-s) &= -\frac{1}{2} + \frac{\log 2 \pi}{2}s + O(s^2) ;\\
\frac{1}{1+s} &= 1-s + O(s^2) ;\\
\lambda^{-s} &=1 -s \log \lambda + O(s^2) ;\\
\zeta(1+s) &= \frac{1}{s} + \gamma + O(s),
\end{align*}
donc le coefficient constant de
$\lambda^{-s}\zeta(-s)\zeta(s+1)/(s+1)$ en $s=0$ est
$$
- \frac{\gamma}{2} + \frac{1}{2} \log \lambda +\frac{1}{2}+ \frac{\log
  2 \pi}{2},
$$
d'o\`u le r\'esultat. \fin

\begin{prop}\label{t129}
Soit $\lambda >0$ tel que la s\'erie $\fhi_1(\lambda)$
converge. Alors,
$$
\fhi_1(\lambda) + \lambda \fhi_1(1/\lambda) -\lambda
\int_{\lambda}^{+\infty}\fhi_1(t) \frac{dt}{t^2} = -\frac{1}{2}
\lambda \log \lambda + \frac{\log
  2 \pi - \gamma}{2}\lambda - \frac{1}{2},
$$
et
$$
\int_0^{+\infty} \frac{\fhi_2(\lambda +u)-\fhi_2(\lambda)}{(\lambda
+u)^3} du = \frac{1}{2} \log \lambda - \frac{\log
  2 \pi - \gamma}{2} + \frac{1}{2\lambda } +\fhi_1(1/\lambda) +
\fhi_1(\lambda)/\lambda.
$$
\end{prop}
\dem

Pour la premi\`ere relation, on compare les propositions \ref{t124} et
\ref{t106} :
$$
\frac{1-\lambda}{2} \log \lambda + \frac{\lambda +1}{2}
(\log 2 \pi - \gamma) - \fhi_1 (\lambda) -\lambda\fhi_1 (1/\lambda)=\frac{1}{2} \log \lambda + \frac{1-\gamma + \log 2\pi
}{2}-\lambda \int_{\lambda}^{+\infty}\fhi_1(t) \frac{dt}{t^2}.
$$

D'autre part, 
\begin{align*}
\int_0^{+\infty} \frac{\fhi_2(\lambda +u)-\fhi_2(\lambda)}{(\lambda
+u)^3} du &=-\frac{\fhi_2(\lambda)}{2 \lambda^2} +
\int_{\lambda}^{+\infty} \fhi_2(u) \frac{du}{u^3}\\
&= \int_{\lambda}^{+\infty} \fhi_1(u) \frac{du}{u^2}\\
&= \frac{1}{2} \log \lambda - \frac{\log
  2 \pi - \gamma}{2} + \frac{1}{2\lambda } +\fhi_1(1/\lambda) +
\fhi_1(\lambda)/\lambda.\fine
\end{align*}

\subsection{D\'eveloppement asymptotique de $A(\lambda)$ au voisinage
d'un nombre rationnel}

Dans ce paragraphe, on consid\`ere deux nombres entiers $p$ et $q$
positifs et
premiers entre eux. Rappelons la notation
$$
\Delta_{p,q}(t) := \fhi_2 \left (\frac{p}{q} +t \right ) - \fhi_2 \left (
\frac{p}{q} \right ), \quad t \in \Real.
$$

\begin{prop}\label{t127}
Pour $t > -p/q$, on a
\begin{multline*}
A \left ( \frac{p}{q} + t \right ) - A \left ( \frac{p}{q}  \right ) =
\frac{1}{2} \log (1 +tq/p) - t \left ( \frac{1}{2} \log
\frac{p}{q} - \frac{\log 2 \pi - \gamma}{2} +\frac{q}{2p} +
\frac{\pi}{2p} \bigl ( V(p,q) + V(q,p) \bigr ) \right ) + \\
\frac{\Delta_{p,q}(t)}{2 \left ( \frac{p}{q} + t \right )} +
\left ( \frac{p}{q} + t \right ) \int_0^t \Delta_{p,q}(u) \left ( \frac{p}{q} + u
\right )^{-3} du.
\end{multline*}
\end{prop}
\dem

D'apr\`es la proposition \ref{t106}, on a
$$
A \left ( \frac{p}{q} + t \right ) = \frac{1}{2} \log \left
( \frac{p}{q} + t \right ) + \frac{1-\gamma + \log 2\pi}{2} +
\frac{\Delta_{p,q}(t)}{2 \left ( \frac{p}{q} + t \right )} - \left
( \frac{p}{q} + t \right ) \int_t^{+\infty} \Delta_{p,q}(u) \left ( \frac{p}{q} + u
\right )^{-3} du.
$$
En particulier,
$$
A \left ( \frac{p}{q}  \right ) = \frac{1}{2} \log  \frac{p}{q}  + \frac{1-\gamma + \log 2\pi}{2} - \frac{p}{q} \int_0^{+\infty} \Delta_{p,q}(u) \left ( \frac{p}{q} + u
\right )^{-3} du,
$$
d'o\`u, par soustraction,
$$
A \left ( \frac{p}{q} + t \right ) - A \left ( \frac{p}{q}  \right )
=\frac{1}{2} \log (1 +tq/p) + \frac{\Delta_{p,q}(t)}{2 \left ( \frac{p}{q} + t \right )} +
\left ( \frac{p}{q} + t \right ) \int_0^t \Delta_{p,q}(u) \left ( \frac{p}{q} + u
\right )^{-3} du - t\int_0^{+\infty} \Delta_{p,q}(u) \left ( \frac{p}{q} + u \right )^{-3}
du,
$$
d'o\`u le r\'esultat, puisque, d'apr\`es la proposition \ref{t129}, on
a
\begin{equation*}
\int_0^{+\infty} \Delta_{p,q}(u) \left ( \frac{p}{q} + u \right )^{-3}
du =\frac{1}{2} \log  \frac{p}{q} -\frac{\log 2\pi - \gamma}{2}
+ \frac{q}{2 p} +\frac{\pi}{2p} \bigl ( V(p,q) + V(q,p) \bigr ).\fine
\end{equation*}

\begin{prop}\label{t103}
Pour $0 <| u| \ioe 1/2q$, on a
\begin{align*}
\Delta_{p,q}(u) & = q^{-1} \bigl ( |u| \log |u| + C^{\pm}(p,q) u -qu^2/2 + O(q^4
u^3) \bigr ) \\
& = q^{-1} \bigl ( |u| \log |u| + C^{\pm}(p,q) u + O(q^3u^2) \bigr
),
\end{align*}
o\`u $\pm$ est le signe de $u$ et
$$
C^{\pm}(p,q) := \pi V(p,q) \pm ( 2 \log q + \log 2 \pi - \gamma -1) \ll q \log
2q.
$$
\end{prop}
{\bf D\'emonstration}

D'apr\`es la proposition \ref{t102}, on a
$$
\Delta_{p,q}(u) = \frac{1}{q} |u| \log |u| + \frac{u}{q} C^{\pm}(p,q) -
\frac{u^2}{2} + O(q^3u^3),
$$
d'o\`u la premi\`ere approximation. D'autre part, $q^3u^3$ et $u^2$
sont $O(q^2u^2)$, ce qui donne la deuxi\`eme. Enfin, $V(p,q) \ll q
\log 2q$ (proposition \ref{t101}), donc $C^{\pm}(p,q)$ v\'erifie la m\^eme estimation.\fin

\begin{prop}\label{t107}
Pour $0 <| t| \ioe 1/2q$, on a
$$
\left ( \frac{p}{q} + t \right )^{-1} \Delta_{p,q}(t) = p^{-1}
\Bigl ( |t| \log |t|  + C^{\pm}(p,q) t \, \mp \,   q p^{-1} t^2 \log |t|- \bigl
(qp^{-1}C^{\pm}(p,q) +q/2 \bigr )t^2 \, \pm \, q^2 p^{-2} t^3 \log |t| + O(q^4t^3)
\Bigr ).
$$
o\`u $\pm$ est le signe de $t$.
\end{prop}
\dem

On a :
\begin{align*}
\left ( \frac{p}{q} + t \right )^{-1} \Delta_{p,q}(t) &= p^{-1} (1 + qt/p)^{-1}\bigl ( |t| \log |t| + C^{\pm}(p,q) t -qt^2/2 + O(q^4
t^3) \bigr )\\
&=p^{-1} \bigl (1 -qt/p + q^2p^{-2}t^2 + O(q^3p^{-3}t^3) \bigr )\bigl ( |t| \log |t| + C^{\pm}(p,q) t -qt^2/2 + O(q^4
t^3) \bigr ).
\end{align*}

De plus,
\begin{align*}
\bigl (1 -qt/p + q^2p^{-2}t^2 + O(q^3p^{-3}t^3) \bigr )O(q^4
t^3) &\ll q^4t^3 ;\\
\bigl (-qt/p + q^2p^{-2}t^2 + O(q^3p^{-3}t^3) \bigr )\cdot (-qt^2/2)
&\ll qp^{-1}t \cdot qt^2 \\
& \ll q^4t^3 ;\\
\bigl ( q^2p^{-2}t^2 + O(q^3p^{-3}t^3) \bigr )\cdot C^{\pm}(p,q)t & \ll
q^2p^{-2}t^2 \cdot q (\log 2q) t\\
& \ll q^4t^3 ;\\
O(q^3p^{-3}t^3) \cdot |t| \log |t| & \ll q^3p^{-3}t^3 \cdot q^{-1} \log 2q
\\
& \ll q^4t^3.
\end{align*}

Dans le d\'eveloppement du produit, il y a donc dix termes dont la
contribution totale est $\ll q^4t^3$, et six termes restants qui sont
:
$$
|t| \log |t| \mp qp^{-1} t^2 \log |t| \pm q^2p^{-2} t^3 \log |t| + C^{\pm}(p,q) t -
qp^{-1}C^{\pm}(p,q)t^2 - q t^2/2,
$$
d'o\`u le r\'esultat. \fin

\begin{prop}\label{t104}
Pour $0 <| u| \ioe 1/2q$, on a
$$
\left ( \frac{p}{q} + u \right )^{-3} \Delta_{p,q}(u) = q^2p^{-3}
\bigl ( |u| \log |u|  + C^{\pm}(p,q) u \mp 3 q p^{-1} u^2 \log |u| + O(q^3u^2) \bigr
).
$$
o\`u $\pm$ est le signe de $u$.
\end{prop}
{\bf D\'emonstration}

\begin{align*}
\left ( \frac{p}{q} + u \right )^{-3} \Delta_{p,q}(u) &= q^3p^{-3} (1
+ qu/p)^{-3}\Delta_{p,q}(u)\\
&= q^2p^{-3} \bigl ( 1 - 3 qu/p + O(q^2p^{-2}u^2) \bigr ) \bigl (  |u| \log |u|  + C^{\pm}(p,q) u + O(q^3u^2) \bigr
),
\end{align*}
d'apr\`es la proposition \ref{t103}. 

On a
$$
1 -3 qu/p + O(q^2p^{-2}u^2) \ll 1,
$$
donc
$$
\bigl ( 1 - 3 qu/p + O(q^2p^{-2}u^2) \bigr ) O(q^3u^2) \ll q^3u^2.
$$

D'autre part,
$$
- 3 qu/p + O(q^2p^{-2}u^2) \ll qp^{-1}u,
$$
donc
\begin{align*}
\bigl ( - 3 qu/p + O(q^2p^{-2}u^2) \bigr )C^{\pm}(p,q) u & \ll q^2 \log 2q
\cdot p^{-1} u^2 \\
& \ll q^3u^2.
\end{align*}

Enfin,
\begin{align*}
O(q^2p^{-2}u^2)|u| \log |u|  & \ll q^2p^{-2}u^3\log |u| \\
&= p^{-2}q^{-1} |u| \log |u| \cdot q^3u^2 \\
& \ll p^{-2}q^{-2} \log 2q \cdot q^3u^2 \\
& \ll q^3u^2.
\end{align*}

Dans le d\'eveloppement du produit, il y a donc six termes dont la
contribution totale est $\ll q^3u^2$, et trois termes restants qui
sont :
\begin{equation*}
|u| \log |u|  + C^{\pm}(p,q) u \mp 3 q p^{-1} u^2 \log |u|.\fine
\end{equation*}

\begin{prop}\label{t105}
Pour $0 < |t| \ioe 1/2q$, on a
\begin{align*}
\int_0^t \left ( \frac{p}{q} + u \right )^{-3} \Delta_{p,q}(u) du &=
\frac{1}{2}q^2p^{-3} \left (\pm t^2\log |t| + \bigl ( C^{\pm}(p,q)\mp 1/2 \bigr ) t^2
\mp 2qp^{-1} t^3 \log |t| + O(q^3t^3) \right )\\
&= q^2p^{-3} \left (\pm  \frac{t^2}{2}\log |t| + O(q^2t^2) \right ).
\end{align*}
o\`u $\pm$ est le signe de $t$.
\end{prop}
\dem

En int\'egrant par parties, on trouve
$$
\int_0^t |u| \log |u| \, du =\pm \left (  \frac{t^2}{2}\log |t| -
\frac{t^2}{4}\right ); \quad
\int_0^t u^2 \log |u| \, du = \frac{t^3}{3}\log |t|- \frac{t^3}{9}.
$$

Par cons\'equent, d'apr\`es la proposition \ref{t104},
\begin{align*}
\int_0^t \left ( \frac{p}{q} + u \right )^{-3} \Delta_{p,q}(u) du &
=q^2p^{-3} \left (\pm \frac{t^2}{2}\log |t| \mp \frac{t^2}{4} +
C^{\pm}(p,q)\frac{t^2}{2}\mp 3qp^{-1} \Bigl ( \frac{t^3}{3}\log |t|-
\frac{t^3}{9} \Bigr )  + O(q^3t^3) \right ) \\
&= \frac{1}{2}q^2p^{-3} \left (\pm  t^2\log |t| + \bigl
( C^{\pm}(p,q)\mp 1/2 \bigr ) t^2
\mp 2qp^{-1} t^3 \log |t| + O(q^3t^3) \right ).
\end{align*}

Pour la deuxi\`eme estimation, on observe que
\begin{align*}
\bigl ( C^{\pm}(p,q)\mp 1/2 \bigr ) t^2 &\ll q (\log 2q) t^2\\
&\ll q^2t^2;\\
qp^{-1} t^3 \log |t| &= qp^{-1} \cdot t \log |t| \cdot t^2\\
&\ll qp^{-1}q^{-1} (\log 2q) t^2\\
&\ll q^2t^2,
\end{align*}
et $q^3t^3 \ll q^2t^2$.\fin

\begin{prop}\label{t108}
Pour $0 < |t| \ioe 1/2q$, on a
$$
\left ( \frac{p}{q} + t \right )\int_0^t \left ( \frac{p}{q} + u
\right )^{-3} \Delta_{p,q}(u) du = \frac{1}{2}qp^{-2} \left ( \pm t^2\log |t| + \bigl ( C^{\pm}(p,q)\mp 1/2 \bigr ) t^2
\mp qp^{-1} t^3 \log |t| + O(q^3t^3) \right ).
$$
o\`u $\pm$ est le signe de $t$.
\end{prop}
\dem

D'apr\`es la proposition \ref{t105}, on a
$$
\frac{p}{q}\int_0^t \left ( \frac{p}{q} + u \right )^{-3} \Delta_{p,q}(u) du =
\frac{1}{2}qp^{-2} \left ( \pm t^2\log |t| + \bigl ( C^{\pm}(p,q)\mp 1/2 \bigr ) t^2
\mp 2qp^{-1} t^3 \log |t| + O(q^3t^3) \right ),
$$
et
$$
t\int_0^t \left ( \frac{p}{q} + u \right )^{-3} \Delta_{p,q}(u) du =
q^2p^{-3} \left ( \pm  \frac{t^3}{2}\log |t| + O(q^2t^3) \right ),
$$
d'o\`u
$$
\left ( \frac{p}{q} + t \right )\int_0^t \left ( \frac{p}{q} + u
\right )^{-3} \Delta_{p,q}(u) du =\frac{1}{2}qp^{-2} \left (\pm t^2 \log |t|
+ \bigl ( C^{\pm}(p,q)\mp 1/2 \bigr ) t^2 \mp (2qp^{-1} -qp^{-1})t^3 \log |t| +
O(q^3t^3) \right ).
$$
\fin

\begin{prop}\label{t128}  
Pour $0 < |t| \ioe 1/2q$, on a
$$
A \left ( \frac{p}{q} + t \right ) - A \left ( \frac{p}{q}  \right ) =
\frac{|t| \log |t|}{2p} + D^{\pm}(p,q) t - \frac{q(p+q\pm 1)}{4p^2} t^2 + O(q^4p^{-1}t^3),
$$
o\`u
$$
D^{\pm}(p,q) := -\frac{1}{2} \log
\frac{p}{q} + \frac{\log 2 \pi - \gamma}{2} \pm \frac{\log q}{p}  \pm
\frac{\log 2 \pi - \gamma -1}{2p} -\frac{\pi}{2p}V(q,p).
$$
et $\pm$ est le signe de $t$.
\end{prop}
\dem

D'apr\`es les propositions \ref{t127}, \ref{t107} et \ref{t108}, on a
\begin{align*}
A \left ( \frac{p}{q} + t \right ) -& A \left ( \frac{p}{q}  \right )
= \frac{1}{2} qp^{-1}t - \frac{1}{4} q^2p^{-2}t^2 + O ( q^3p^{-3}t^3)
+\\
& \quad - t \left ( \frac{1}{2} \log
\frac{p}{q} - \frac{\log 2 \pi - \gamma}{2} +\frac{q}{2p}+
\frac{\pi}{2p} \bigl ( V(p,q) + V(q,p) \bigr ) \right ) + \\
& \quad + \frac{1}{2}p^{-1}
\Bigl ( |t| \log |t|  + C^{\pm}(p,q) t \mp  q p^{-1} t^2 \log |t|- \bigl
(qp^{-1}C^{\pm}(p,q) +q/2 \bigr )t^2 \pm q^2 p^{-2} t^3 \log |t| + O(q^4t^3)
\Bigr ) +\\
& \quad  \frac{1}{2}qp^{-2} \left ( \pm t^2\log |t| + \bigl ( C^{\pm}(p,q)\mp 1/2 \bigr ) t^2
\mp qp^{-1} t^3 \log |t| + O(q^3t^3) \right ).
\end{align*}

Le coefficient de $|t| \log |t|$ est $1/2p$; ceux de $t^2 \log |t|$ et $t^3
\log |t|$ sont nuls. Le coefficient de t est :
\begin{align*}
& \quad \frac{1}{2} qp^{-1}-  \left ( \frac{1}{2} \log
\frac{p}{q} - \frac{\log 2 \pi - \gamma}{2} +\frac{q}{2p} +
\frac{\pi}{2p} \bigl ( V(p,q) + V(q,p) \bigr ) \right ) +
\frac{1}{2}p^{-1}C^{\pm}(p,q)\\
&= -\frac{1}{2} \log
\frac{p}{q} + \frac{\log 2 \pi - \gamma}{2} \pm \frac{\log q}{p}  \pm
\frac{\log 2 \pi - \gamma -1}{2p} -\frac{\pi}{2p}V(q,p).
\end{align*}

Le coefficient de $t^2$ est
\begin{align*}
& \quad - \frac{1}{4} q^2p^{-2}  - \frac{1}{2}p^{-1} \bigl (qp^{-1}C^{\pm}(p,q) +q/2
  \bigr ) + \frac{1}{2}  qp^{-2}\bigl ( C^{\pm}(p,q)\mp 1/2 \bigr )\\
&= -\frac{1}{4} qp^{-2} (p+q\pm 1).
\end{align*}

Enfin, le terme compl\'ementaire est $O(q^4p^{-1} t^3)$.\fin

\end{document}